\pdfoutput=1
\documentclass[11pt]{article}
\usepackage[margin=1in]{geometry}
\usepackage{amsmath,amssymb,amsthm,mathtools}
\usepackage{booktabs,tabularx,array}
\usepackage{float}
\usepackage{placeins}
\usepackage{enumitem}
\usepackage{xcolor}
\usepackage[numbers,sort&compress]{natbib}
\usepackage[hypertexnames=false]{hyperref}
\usepackage{microtype}
\usepackage{graphicx}
\usepackage{wrapfig}
\usepackage{capt-of}
\graphicspath{{figures/}}
\usepackage[capitalize,nameinlink,noabbrev]{cleveref}
\crefname{assumption}{Assumption}{Assumptions}
\Crefname{assumption}{Assumption}{Assumptions}
\hypersetup{colorlinks=true,linkcolor=blue!55!black,citecolor=blue!55!black,urlcolor=blue!55!black}
\allowdisplaybreaks
\numberwithin{equation}{section}
\newtheorem{theorem}{Theorem}[section]
\newtheorem{lemma}[theorem]{Lemma}
\newtheorem{proposition}[theorem]{Proposition}
\newtheorem{corollary}[theorem]{Corollary}
\theoremstyle{definition}
\newtheorem{assumption}{Assumption}
\newtheorem{definition}[theorem]{Definition}
\newtheorem{example}[theorem]{Example}
\theoremstyle{remark}
\newtheorem{remark}[theorem]{Remark}

\newcommand{\by}[2]{\overset{\text{\scriptsize #1}}{#2}}
\newsavebox{\estimatorcontent}
\newenvironment{estimatorbox}[1]{%
  \par\medskip\noindent
  \begin{lrbox}{\estimatorcontent}%
  \begin{minipage}{\dimexpr\linewidth-2\fboxsep-2\fboxrule\relax}
  \textbf{#1.}\par\smallskip
}{%
  \end{minipage}\end{lrbox}%
  \fbox{\usebox{\estimatorcontent}}\par\medskip
}
\title{First-Order Information from Noisy Comparison Oracles}
\author{%
Qipei Chen\thanks{Authors in alphabetical order.}%
\textsuperscript{\, }\thanks{\raggedright Department of Mathematics, University of
California, Berkeley.
Email: \href{mailto:perrychen@berkeley.edu}{perrychen@berkeley.edu}.
Chen contributed to the earlier stages of this work before summer 2026.
The vast majority of his contributions were made while he was at UC Berkeley.}
\qquad
Katya Scheinberg\thanks{\raggedright H. Milton Stewart School of Industrial and Systems
Engineering, Georgia Institute of Technology.
Email: \href{mailto:katya.scheinberg@isye.gatech.edu}{katya.scheinberg@isye.gatech.edu}.}
\qquad
Zikai Xiong\thanks{\raggedright Department of Industrial Engineering and Management Sciences,
Northwestern University.
Email: \href{mailto:zikai.xiong@northwestern.edu}{zikai.xiong@northwestern.edu}
(corresponding author).}}
\date{November 22, 2026}
\begin{document}
\maketitle

\begin{abstract}
We study stochastic and deterministic smooth unconstrained optimization problems
where neither the functions nor their derivatives can be computed or estimated and the only information
about the objective is available via pairwise comparisons.
We focus on a specific model of smooth comparison oracles, as in the
Bradley--Terry \citep{bradley1952rank} and Thurstone \citep{thurstone1927law} models, in which the probability of a correct
comparison depends smoothly on the function-value gap.
We investigate what first-order information can be obtained from
such comparisons in deterministic and stochastic settings, and how many
comparisons are needed to estimate it. We show that the
first-order information that such comparisons provide is a fixed positive
multiple of the gradient of the objective. We propose a random-direction
estimator of the comparison gradient and analyze its bias and variance. It
reaches mean-squared error $\varepsilon^2$ with $O(d^2\varepsilon^{-4})$
comparisons when the sample functions have Lipschitz gradients, and
$O(d^2\varepsilon^{-3})$ when their Hessians are also Lipschitz. We prove
lower bounds of the same order in $d$ and $\varepsilon$ for all
comparison-based estimators, so these numbers are necessary and our estimator
is optimal in $d$ and $\varepsilon$. As an application, we use it in a
stochastic gradient method to find a point that is approximately stationary
with respect to comparisons. We also extend the estimator to nonsmooth
objectives, for which it estimates a smoothed comparison gradient.
\end{abstract}

\section{Introduction}\label{sec:intro}

We study the minimization of a function $f\colon\mathbb R^d\to\mathbb R$
whose values cannot be observed. The only available information comes from
pairwise comparisons that report, possibly incorrectly, which of two given points has
the smaller function value. We ask what first-order information at a given
point can be obtained from such comparisons and how many comparisons are
needed to estimate it.

Comparisons are widely used when the quality of a decision is assessed by
people, as in reinforcement learning from human preferences
\citep{Christiano2017HumanPreferences} and in the fine-tuning of language
models \citep{Ouyang2022InstructGPT,Rafailov2023DPO}. In a variety of other applications involving human
feedback an answer that compares two decisions is often more reliable than
assigning each decision a numerical value (see, e.g., \citet{ShahEtAl2016Topology}). For
example, it is easier for a person to say which of two designs looks better
than to give each of them a score or to evaluate how much better one is than the other numerically \citep{BrochuEtAl2007ActivePreference}. Moreover,
the number of such applications is growing, since language models
can now mimic human behavior, including making comparisons, at a much smaller
cost \citep{Zheng2023LLMJudge,Dubois2023AlpacaFarm,Lee2024RLAIF}.

In many of these applications, the outcome of a decision also depends on a
random sample. For example, the outputs of two search engines with different internal parameters
are compared through their results for the (random) same query \citep{yue2009interactively}.
Similarly, two language models can be compared through their responses to the
same (random) prompt \citep{Munos2024NLHF}.
In all of these applications, one looks for the best decision (e.g. parameters defining the language model),
which results in optimization problems whose objective function can only be observed through
comparisons. Such problems are stochastic when the objective function is an expectation over
random samples, in which case each comparison is made on a single sample.

Thus, in this paper, we study the unconstrained optimization problem
\begin{equation}
 \min_{x\in\mathbb R^d} f(x)\,,
 \tag{D}\label{eq:deterministic-problem}
\end{equation}
and its stochastic variant
\begin{equation}
 \min_{x\in\mathbb R^d} f(x)\,,\qquad
 f(x):=\mathbb E_{\xi\sim\mathbb P_\xi}[f(x;\xi)]\,,
 \tag{S}\label{eq:problem}
\end{equation}
where $\xi$ is a sample from a distribution $\mathbb P_\xi$ on $\Xi$, such as
a random query or a prompt, or a random responder from a population. For each $\xi$ we assume the existence
of the cost function $f(x;\xi)$ which is revealed through pairwise comparisons but cannot be evaluated.
We will assume certain properties of $f(x)$ and $f(x;\xi)$ which we will discuss later.

\begin{figure}[htbp]
\centering
\includegraphics[width=\textwidth]{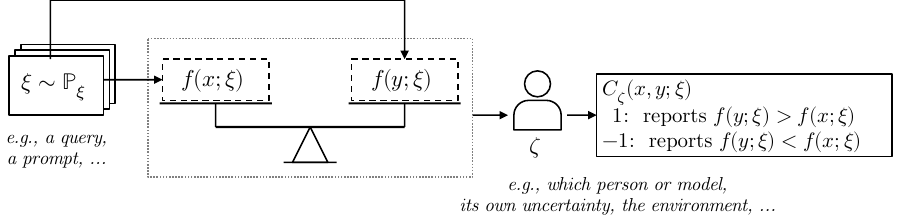}
\caption{One comparison in problem~\eqref{eq:problem}.}
\label{fig:one-comparison}
\end{figure}

Let us now introduce our key notation and concepts.
A \emph{comparison oracle} queried at two points $x$ and
$y$ returns a possible random outcome $C_\zeta(x,y)\in\{-1,1\}$, where $C_\zeta(x,y)=1$ is understood as
$f(y)>f(x)$. The output of the oracle
is written as $C_\zeta(x,y;\xi)$, when we wish to indicate that the comparison is
performed for $f(x;\xi)$ and $f(y;\xi)$ for a particular sample $\xi$. The random variable
$\zeta$ represents the internal randomness of comparison oracle which can come from the
 uncertainty or from the environment, which may cause the output to be incorrect.
 For example if a user is asked to compare two very similar designs they may give a random
 response depending on which design they see first or some other factor unrelated to design "quality".
 The probability of an oracle output error
 will be a result of the realization of $\zeta$ and may depend on the function $f(x)$ or $f(x;\xi)$.
\Cref{fig:one-comparison} shows the flow of the comparison oracle in
problem~\eqref{eq:problem} indicating how the two sources of randomness, $\xi$ and
$\zeta$ are realized.

\begin{figure}[htbp]
\centering
\includegraphics[width=\textwidth]{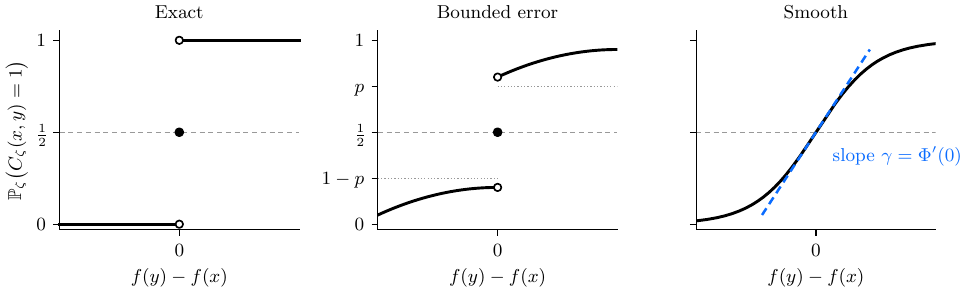}
\caption{The probability $\mathbb P_\zeta(C_\zeta(x,y)=1)$ under the three models
in~\eqref{eq:oracle-models}.}
\label{fig:three-models}
\end{figure}

We now define three standard models of comparison oracles and illustrate them in
\cref{fig:three-models}. Each model specifies the probability that the oracle output is correct. In particular, when
$f(x)\ne f(y)$ the oracle satisfies one of the three conditions stated as follows
\begin{equation}
 \mathbb P_\zeta \bigl(C_\zeta (x,y)=\operatorname{sign}(f(y)-f(x))\bigr)
 \begin{cases}
  =1 & \text{(exact)},\\
  \ge p\ \text{ for a fixed } p>\tfrac12 & \text{(bounded error)},\\
  =\tfrac12+\Phi\bigl(|f(y)-f(x)|\bigr) & \text{(smooth)},
 \end{cases}
\label{eq:oracle-models}
\end{equation}
with $f$ replaced by $f(\cdot;\xi)$ and $C_\zeta (x,y)$ replaced by $C_\zeta (x,y;\xi)$ for problem~\eqref{eq:problem}.

 The bounded-error comparison model goes back to the noisy comparisons studied in
\citet{feige1994noisy}. The comparison oracle under such model is correct with probability at least $p$ regardless
of how close the two compared vectors are in terms of the objective value.
Exact comparisons have no errors and are a special case of bounded-error
comparisons with $p=1$. On the other hand, bounded-error comparisons can be made nearly exact
by applying the same comparison oracle repeatedly and choosing the majority response. Such an oracle will be correct with overwhelmingly high probability.

Exact comparisons have long been used in
optimization. For example, direct-search methods in derivative-free optimization, such as those of
\citet{HookeJeeves1961} and \citet{NelderMead1965}, use only comparisons of
function values (see, e.g., \citet{KoldaLewisTorczon2003DirectSearch} and
\citet{ConnScheinbergVicente2009}), and \citet{jamieson2012query} formalized
comparisons as an oracle. More recently, exact comparison oracles have been
used in optimization of parameters of a language model
\citep{ChenEtAl2025ComPO,ChenEtAl2026ZerothOrderAlignment}. Exact
comparisons depend only on the ordering of function values, and in
problem~\eqref{eq:deterministic-problem} this ordering is unchanged when $f$
is multiplied by a positive constant. This implies that exact comparisons alone cannot identify
the function $f$ or its gradient $\nabla f(x)$. Nevertheless, it is possible to
approximate the gradient direction $n(x):=\tfrac{\nabla f(x)}{\|\nabla f(x)\|}$ using only
exact comparisons, and \citet{tao2026gradient} and
\citet{scheinberg2026functionfree} estimate $n(x)$ with a number of
comparisons that grows only logarithmically in the required accuracy.

Exact and nearly-exact comparison oracles do not capture many realistic applications.
For example, when a person or a language model is comparing two alternatives (generated by $x$ and $y$)
such that their objective values $f(x;\xi)$ and $f(y;\xi)$ are very close, these two alternatives may be nearly indistinguishable
and thus are more likely to be compared incorrectly. The smooth comparison model
in~\eqref{eq:oracle-models} captures this behavior via the use of function $\Phi$, which is a smooth
odd function with $\Phi'(0)=\gamma>0$, such that
$\mathbb P_\zeta(C_\zeta(x,y)=1)=\tfrac12+\Phi(f(y)-f(x))$. Through $\Phi$, the probability of a correct
comparison depends on the size of the function-value gap. In
particular, because $\Phi(0)=0$, a comparison of two nearly equal values is
close to a fair coin flip. This model includes the Thurstone model
\citep{thurstone1927law} and the Bradley--Terry model
\citep{bradley1952rank} of paired comparisons. The Bradley--Terry model has
been widely used to model human preferences when language models are
fine-tuned from comparisons \citep{Ouyang2022InstructGPT,Rafailov2023DPO}.

The key mathematical difference of the smooth comparison oracle model from the exact and bounded-error models
is that the probability of the correct outcome of the former carries some information about the function gap. This information is encoded
via $\Phi$. Past works that considered this model make different assumptions about the level of knowledge of $\Phi$.

Smooth comparisons have been used to minimize
deterministic convex objectives, with bounds on the regret or on the error in
function value in
\citep{yue2009interactively,kumagai2017regret,saha2025general} without assuming
the knowledge of $\Phi$ itself but assuming a bound on its growth.

In \citet{LuEtAl2026DuelingFeedback} it is assumed that $\Phi'(0)$ is known, while in
\citet{ZhangYing2025ZPG} and \citet{ChaytiElBakkaliSaadi2026NoisyComparisons}
it is assumed that $\Phi$ itself is known. Using the analytical form of $\Phi$ it is
possible to recover the differences of function values
from comparison oracles through the inverse of $t\mapsto\tfrac12+\Phi(t)$. Similarly,
one can recover the differences of function values locally using the explicit knowledge of $\Phi'(0)$.

However, none of these works optimize a stochastic objective for a general smooth $\Phi$.
Moreover, in all of these works the goal is stated in terms of $f$, such as
its regret, its function values, or its stationary points. Comparisons,
however, without knowledge of $\Phi$ do not determine $f$ and $\Phi$ separately, because multiplying $f$
by a positive constant and rescaling $\Phi$ accordingly leaves the
probability of every comparison unchanged. Since assuming that $\Phi$ is known and fixed
is too strong and limiting, we claim that the goal of optimization based on smooth
comparison oracles should not be stated in terms of $f$ and should be stated via information
that is pertinent to comparison oracles.

This leads to the following two questions.
\begin{enumerate}
 \item[Q1.] \textit{For each model in~\eqref{eq:oracle-models}, what
 first-order information relevant to
 the deterministic problem~\eqref{eq:deterministic-problem} or the stochastic
 problem~\eqref{eq:problem}, at a given $x$, can be identified using only comparisons?}
 \item[Q2.] \textit{For the information that is identifiable, how many
 comparisons are necessary and sufficient to estimate it?}
\end{enumerate}
For exact and bounded-error comparisons in
problem~\eqref{eq:deterministic-problem}, both questions are answered by the
works of \citet{tao2026gradient} and \citet{scheinberg2026functionfree}
discussed above. For problem~\eqref{eq:problem}, we show that such comparisons
may not identify even the gradient direction. For smooth comparisons, this
paper answers both questions for both problems.

\subsection{Outline and contributions}
\label{subsec:outline-contributions}

In \cref{sec:comparison-information}, we answer Q1. Under smooth
comparisons, the probability that $x$ is preferred to a nearby point $y$,
which we call the \emph{comparison probability}, changes smoothly with $y$,
and the first-order information that comparisons provide is its gradient
\[
 G(x):=\nabla_y\,\mathbb P_\zeta\bigl(C_\zeta(x,y)=1\bigr)\big|_{y=x}\,,
\]
which we call the \emph{comparison gradient}. In problem~\eqref{eq:problem},
$G(x)$ is the average of the comparison gradients of the sample functions. In both
deterministic and stochastic problems, it equals $\gamma\nabla f(x)$ with
$\gamma=\Phi'(0)$, but $\nabla f(x)$ itself is not identifiable. The comparison gradient describes the rate at
which the probability that a nearby point is preferred to $x$ moves away from
$\tfrac12$. When it is small, no small move away from $x$ is preferred with
probability much above $\tfrac12$, and in this sense $x$ is nearly stationary
with respect to comparisons.

In \cref{sec:estimation}, we estimate $G(x)$ up to mean-squared error $\varepsilon^2$.
We use comparisons of $x$ and $x+hu$, with random directions $u$, which is a technique borrowed from
well-studied finite-difference estimators based on smoothing techniques \citep{NesterovSpokoiny2017,FlaxmanKalaiMcMahan2005,BerahasCaoChoromanskiScheinberg2022}.
Neither our method nor our analysis assumes knowledge of
$\Phi$ or $\gamma$. In both deterministic and stochastic optimization, we show that
our estimator reaches this accuracy with at most
$O(d^2\varepsilon^{-4})$ comparisons when $\nabla f$ is Lipschitz, and with
$O(d^2\varepsilon^{-3})$ comparisons when $\nabla^2 f$ is also Lipschitz. In
\cref{sec:gradient-lower-bound}, we show how many comparisons are necessary for this task by
proving lower bounds of the same order in $d$ and $\varepsilon$ that hold
even for comparison-based estimators that know $\Phi$. These lower bounds
describe a limit of the information contained in comparisons rather than of a
particular method, and since our proposed estimator attains them, its dependence on
$d$ and $\varepsilon$ cannot be improved. \Cref{tab:results}
summarizes the answers to Q1 and Q2.

In \cref{sec:optimization-application}, we use the estimator in a standard
stochastic gradient framework to find a point that is nearly stationary with
respect to comparisons, with $O(d^2\varepsilon^{-6})$ comparisons in
expectation when the sample functions have Lipschitz gradients and
$O(d^2\varepsilon^{-5})$ when their Hessians are also Lipschitz. In
\cref{sec:nonsmooth-estimation}, we consider sample functions that are
Lipschitz but not necessarily differentiable, and we estimate the smoothed
comparison gradient (which is $\gamma$ times the gradient of a Gaussian
smoothing of $f$ with parameter $h>0$) to mean-squared error $\varepsilon^2$
with $O\bigl(d^2\varepsilon^{-2}\max\{h^{-2},\varepsilon^{-1}\}\bigr)$
comparisons.

\begin{table}[htbp]
\centering
\small
\setlength{\abovecaptionskip}{0pt}
\setlength{\belowcaptionskip}{6pt}
\caption{Numbers of comparisons sufficient and necessary for estimating
first-order information.}
\label{tab:results}
\begin{tabular*}{\textwidth}{@{\extracolsep{\fill}}llll@{}}
\toprule
Oracle & Identified & Sufficient & Necessary \\
\midrule
\multicolumn{4}{@{}l}{\textit{Deterministic problem}~\eqref{eq:deterministic-problem}} \\
Exact or bounded-error
 & $n(x)$
 & $O(d\log(\tfrac d\varepsilon))$~\citep{scheinberg2026functionfree,tao2026gradient}
 & $\Omega(d\log(\tfrac1\varepsilon))$~\citep{scheinberg2026functionfree,tao2026gradient} \\
Smooth
 & $G(x)$
 & $O(d^2\varepsilon^{-4})$, $O(d^2\varepsilon^{-3})^{\dagger}$
 & $\Omega(d^2\varepsilon^{-4})$, $\Omega(d^2\varepsilon^{-3})^{\dagger}$ \\
\addlinespace
\multicolumn{4}{@{}l}{\textit{Stochastic problem}~\eqref{eq:problem}} \\
Exact or bounded-error
 & Neither $\nabla f(x)$ nor $n(x)$
 & --
 & -- \\
Smooth
 & $G(x)$
 & $O(d^2\varepsilon^{-4})$, $O(d^2\varepsilon^{-3})^{\dagger}$
 & $\Omega(d^2\varepsilon^{-4})$, $\Omega(d^2\varepsilon^{-3})^{\dagger}$ \\
\bottomrule
\end{tabular*}

\smallskip
{\footnotesize\raggedright
$^{\dagger}$The second rate assumes, in addition, that the Hessian is Lipschitz.\par}
\end{table}

\subsection{Notation}
\label{subsec:notation}

For vectors, $\|\cdot\|$ denotes the Euclidean norm and
$\langle\cdot,\cdot\rangle$ the Euclidean inner product. We use $I_d$ to
denote the $d\times d$ identity matrix, $\operatorname{tr}(A)$ the trace of a
matrix $A$, and $\mathbb B^d$ and $\mathbb S^{d-1}$ the closed unit ball and
the unit sphere in $\mathbb R^d$. For a real-valued function $g$, $\nabla g$
and $\nabla^2 g$ denote its gradient and Hessian, and a subscript on $\nabla$
or $\partial$ (e.g., $\nabla_y$ or $\partial_\theta$) indicates the variable of
differentiation.

We use $N(\mu,\Sigma)$ to denote the Gaussian distribution with mean $\mu$ and
covariance $\Sigma$. A subscript on $\mathbb E$ or $\mathbb P$ indicates the
random variables over which the expectation or probability is taken, and
$\operatorname{sign}(t)\in\{-1,0,1\}$ is the sign of $t\in\mathbb R$. A
reference written above an equality or inequality indicates the result that
justifies it. When we state the dependence of a number of comparisons on $d$
and $\varepsilon$, the notation $O(\cdot)$ and $\Omega(\cdot)$ hides factors
that do not depend on $d$ or $\varepsilon$.

\section{Comparison oracles and first-order information}
\label{sec:comparison-information}
In this section, we state the three models in~\eqref{eq:oracle-models}
with the conditions used in the analysis and determine what first-order
information each of them identifies
for the problems~\eqref{eq:deterministic-problem} and~\eqref{eq:problem}.

\subsection{Comparison oracles}
\begin{definition}[Comparison oracle]\label[definition]{def:comparison-oracle}
Let $\zeta$ be a random variable sampled independently at each oracle call.
A comparison oracle queried at $x,y\in\mathbb R^d$ returns
$C_\zeta(x,y)\in\{-1,1\}$ that is determined by $x,y,\zeta$, with
$C_\zeta(y,x)=-C_\zeta(x,y)$ for $x\ne y$ and each realization of $\zeta$.
When $f(x)\ne f(y)$, the comparison is \emph{correct} if
\[
 C_\zeta(x,y)=\operatorname{sign}(f(y)-f(x))\,.
\]
\end{definition}
As in \cref{fig:one-comparison}, $\zeta$ represents the internal randomness of
the oracle. For problem~\eqref{eq:problem}, the oracle compares $x$ and $y$ on
a sample function $f(\cdot;\xi)$ and returns $C_\zeta(x,y;\xi)$, which is
determined by $x,y,\xi,\zeta$. The definitions and conditions in this
section are stated for problem~\eqref{eq:deterministic-problem}, and for
problem~\eqref{eq:problem} they apply to each sample function, with
$f(\cdot;\xi)$ and $C_\zeta(x,y;\xi)$ in place of $f$ and $C_\zeta(x,y)$ and
with probabilities conditional on $\xi$. The comparison probability is
$\mathbb P_\zeta(C_\zeta(x,y)=1\mid x,y)$, which we write as
$\mathbb P(C(x,y)=1)$. In problem~\eqref{eq:problem}, the
comparison probability of the sample function $f(\cdot;\xi)$ is
$\mathbb P(C(x,y;\xi)=1):=\mathbb P_\zeta(C_\zeta(x,y;\xi)=1\mid x,y,\xi)$,
which averages over $\zeta$ but not over $\xi$. We abbreviate other
probabilities of responses in the same way, such as the probability $\mathbb P(C(x,y)=\operatorname{sign}(f(y)-f(x)))$ of a
correct comparison.

\subsection{Exact and bounded-error comparisons}
\label{subsec:exact-bounded-error}
\cref{def:comparison-oracle} does not specify the probability of a correct
comparison. A classical condition of \citet{jamieson2012query} bounds this
probability in terms of the function-value gap: for $f(x)\ne f(y)$,
\begin{equation}
 \mathbb P\bigl(C(x,y)=\operatorname{sign}(f(y)-f(x))\bigr)
 \ge\frac12+\min\{\delta_0,\mu|f(y)-f(x)|^\beta\}\,,
 \label{eq:bounded-error}
\end{equation}
where $\delta_0\in(0,\tfrac12]$, $\mu>0$, and $\beta\ge0$. When $\beta=0$, the
right-hand side equals $\tfrac12+\min\{\delta_0,\mu\}$ for all $x$ and $y$, so
\eqref{eq:bounded-error} reduces to the bounded-error
condition defined below. Although
\eqref{eq:bounded-error} is more general,
comparison-based methods are often analyzed only for exact or bounded-error
comparisons \citep{Bergou2020STP,karabag2021smooth,saha2021dueling}. For
example, \citet{cai2022onebit} state their oracle
under~\eqref{eq:bounded-error}, but their analysis uses
only the case $\beta=0$.

\begin{definition}[Exact and bounded-error comparison oracles]
\label[definition]{def:exact-bounded-error}
A \emph{bounded-error comparison oracle} is a comparison oracle from \cref{def:comparison-oracle}
for which there exists $p\in(\tfrac12,1]$ such that, for any $x,y$ with
$f(x)\ne f(y)$,
\begin{equation*}
 \mathbb P\bigl(C(x,y)=\operatorname{sign}(f(y)-f(x))\bigr)\ge p\,.
\end{equation*}
The oracle is called an \emph{exact comparison oracle} if this
condition holds with $p=1$.
If the queried function values are equal, each response in $\{-1,1\}$
has probability $\tfrac12$.
\end{definition}
A bounded-error comparison oracle can be reduced to an exact one with high
probability. Repeating a bounded-error comparison
of the same pair and taking a majority vote gives the correct ordering with
probability at least $1-\delta$ after
$O\bigl((p-\tfrac12)^{-2}\log((p-\tfrac12)^{-2}\delta^{-1})\bigr)$
repetitions \citep{saha2021dueling}. This also holds in the stochastic
problem when the repetitions use the same sample function. We therefore
discuss exact and bounded-error comparisons together.

Exact and bounded-error comparisons have also been
used to estimate the gradient direction $n(x)$
\citep{karabag2021smooth,cai2022onebit}.
It was shown in \citet{tao2026gradient} and \citet{scheinberg2026functionfree}
that for $\varepsilon\le\tfrac1d$, order $d\log(\tfrac1\varepsilon)$
exact comparisons are sufficient and necessary to compute $\hat n$ satisfying $\|\hat n-n(x)\|\le\varepsilon$.

We now show that $\nabla f(x)$ itself is not identifiable in
problem~\eqref{eq:deterministic-problem}. Since exact oracles satisfy
\cref{def:exact-bounded-error}
and~\eqref{eq:bounded-error} for all values of the
parameters, it suffices to show this under exact comparisons.

\begin{proposition}
\label[proposition]{prop:deterministic-nonidentifiability}
Under exact comparison oracles, any gradient estimator at zero has
mean-squared error at least $\tfrac14$ for one of the functions $f_1(z):=z$
and $f_2(z):=2z$ on $\mathbb R$.
\end{proposition}
\begin{proof}
Since $f_2=2f_1$, the two functions order any two points in the same way, so,
under exact comparison oracles, every comparison returns $1$ with the same
probability for $f_1$ and for $f_2$. Since each query depends only on the
estimator's own random choices and the previous responses, the queries, the
responses, and $\hat g$ have the same distribution for $f_1$ and for $f_2$.
As $\nabla f_1(0)=1$, $\nabla f_2(0)=2$, and $(g-1)^2+(g-2)^2\ge\tfrac12$ for
every $g\in\mathbb R$, one of the two mean-squared errors is at least
$\tfrac14$.
\end{proof}

In stochastic optimization, both exact and bounded-error comparisons
 depend on each sample function only through its
ordering of points. Such comparisons may not even be useful to identify the gradient direction information
for the expectation. This is
because rescaling the individual sample functions by different positive
factors preserves their orderings but can change the gradient direction of
their expectation. \citet{saha2021dueling} use this observation to show that
the minimizer of the objective may not be identifiable. The following
proposition shows that the gradient direction may not be identifiable either.
\begin{proposition}
\label[proposition]{prop:stochastic-nonidentifiability}
Let $\xi$ be uniform on $\{1,2\}$, and consider the two stochastic problems with
sample functions
\begin{equation*}
 f_1(z;1):=3z\,,\quad f_1(z;2):=-z
 \qquad\text{and}\qquad
 f_2(z;1):=z\,,\quad f_2(z;2):=-3z\,,
 \qquad z\in\mathbb R\,.
\end{equation*}
Their objectives $f_1(z)=z$ and $f_2(z)=-z$ have gradients $1$ and $-1$ at
zero. Under exact comparison oracles, any gradient estimator at zero has
mean-squared error at least $1$ for one of the two problems, and any estimator
$\hat n\in\{-1,1\}$ of the gradient direction at zero takes the wrong sign with
probability at least $\tfrac12$ for one of them.
\end{proposition}
\begin{proof}
For each $\xi$, the two sample functions are positive multiples of each other,
so, as in the proof of \cref{prop:deterministic-nonidentifiability}, $\hat g$ and $\hat n$
have the same distribution in the two problems. As
$(g-1)^2+(g+1)^2\ge2$ for every $g\in\mathbb R$, one of the two mean-squared
errors is at least $1$, and the two error probabilities of $\hat n$ sum to
$\mathbb P(\hat n=-1)+\mathbb P(\hat n=1)=1$.
\end{proof}
By \cref{prop:stochastic-nonidentifiability}, in stochastic optimization,
condition~\eqref{eq:bounded-error} does not guarantee that $\nabla f(x)$ or
even $n(x)$ can be determined, however many comparisons are used, and neither
does any other condition that only bounds the probability of a correct
comparison from below, since exact oracles satisfy all such conditions.

\subsection{Smooth comparisons and the comparison gradient}
\label{subsec:smooth-comparisons}
We now define smooth comparisons in~\eqref{eq:oracle-models} with the
specific conditions on $\Phi$ used in the analysis.
\begin{definition}[Smooth comparison oracle]\label[definition]{def:smooth-oracle}
A comparison oracle in \cref{def:comparison-oracle} is called a
\emph{smooth comparison oracle with parameters $(\gamma,M_\phi)$},
where $\gamma>0$ and $M_\phi>0$, if there is an odd, three times
differentiable function $\Phi:\mathbb R\to[-\tfrac12,\tfrac12]$ such that
\begin{equation}
 \mathbb P(C(x,y)=1)
 =\frac12+\Phi(f(y)-f(x))\,.
 \label{eq:comparison-probability}
\end{equation}
The function $\Phi$ satisfies
\begin{equation*}
 \Phi'(0)=\gamma\,,\qquad\text{and}\qquad |\Phi'''(t)|\le M_\phi
 \quad\text{for all }t\in\mathbb R\,.
\end{equation*}
\end{definition}
Note that $\Phi(0)=\Phi''(0)=0$ because $\Phi$ is odd, so each
response has probability $\tfrac12$ when the queried function values are
equal. Moreover, for any $t\in\mathbb R$ it holds that
\begin{equation}
 |\Phi''(t)|\le M_\phi|t|\,,
 \qquad |\Phi'(t)-\gamma|\le\frac{M_\phi}{2}|t|^2\,,
 \qquad\text{and}\qquad |\Phi'(t)|\le\frac23M_\phi^{1/3}\,.
\label{eq:Phi-bounds}
\end{equation}
These bounds are used in the proofs below, and we prove them in
\hyperref[app:proofs]{the appendix}. In particular, $t=0$ in the last bound
gives $M_\phi\ge\tfrac{27}{8}\gamma^3$. In problem~\eqref{eq:problem}, $\Phi$
may depend on $\xi$, while $\gamma$ and $M_\phi$ are the same for all $\xi$.

The two classical models of paired comparisons, the Bradley--Terry and
Thurstone models, are smooth comparison oracles, as the following two examples
show.

\par\medskip\noindent
\begin{minipage}[c]{0.5\textwidth}
\begin{example}[Bradley--Terry comparisons]\label[example]{ex:bradley-terry}
In the Bradley--Terry model \citep{bradley1952rank},
$\mathbb P(C(x,y)=1)=\tfrac{1}{1+\exp(-4\gamma(f(y)-f(x)))}$. This is
\eqref{eq:comparison-probability} with $\Phi(t)=\tfrac12\tanh(2\gamma t)$,
which satisfies \cref{def:smooth-oracle} with $M_\phi=8\gamma^3$.
\end{example}

\begin{example}[Thurstone comparisons]\label[example]{ex:thurstone}
In the Thurstone model \citep{thurstone1927law}, the oracle compares
$f(x)+\sigma Z_x$ and $f(y)+\sigma Z_y$, where $Z_x$ and $Z_y$ are independent
standard normal variables drawn at each call. Then
$\mathbb P(C(x,y)=1)
=\mathbb P\bigl(Z\le\tfrac{f(y)-f(x)}{\sqrt2\sigma}\bigr)$ with $Z\sim N(0,1)$, which
yields $\Phi(t)=\mathbb P(Z\le\tfrac{t}{\sqrt2\sigma})-\tfrac12$. This satisfies
\cref{def:smooth-oracle} with $\gamma=\tfrac{1}{2\sqrt\pi\,\sigma}$ and
$M_\phi=2\pi\gamma^3$.
\end{example}
\end{minipage}\hfill
\begin{minipage}[c]{0.47\textwidth}
\centering
\includegraphics[width=\linewidth]{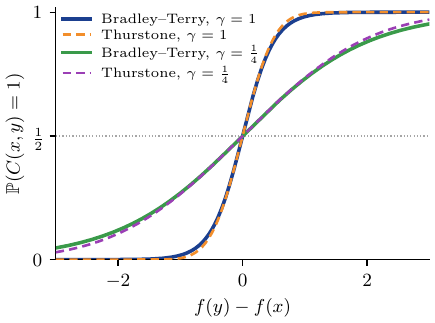}
\captionof{figure}{The Bradley--Terry and Thurstone models for two values of $\gamma$.}
\label{fig:bt-thurstone}
\end{minipage}
\par\vspace{\intextsep}

Both models are widely used, the Bradley--Terry model in ranking from paired
comparisons \citep{Hunter2004MM} and in fine-tuning language models from human
preferences \citep{Ouyang2022InstructGPT,Rafailov2023DPO}, and the Thurstone
model in skill rating \citep{Herbrich2006TrueSkill} and preference learning
\citep{ChuGhahramani2005Preference}.

By~\eqref{eq:comparison-probability}, the comparison probability depends
smoothly on $y$ when $f$ does, and its gradient at $y=x$ defines the
comparison gradient.

\begin{definition}[Comparison gradient]\label[definition]{def:comparison-gradient}
The \emph{comparison gradient} at $x\in\mathbb R^d$ is
\begin{equation}
 G(x):=\nabla_y\,\mathbb P\bigl(C(x,y)=1\bigr)\big|_{y=x}\,.
 \label{eq:comparison-gradient}
\end{equation}
In problem~\eqref{eq:problem}, it is the average of the comparison gradients
of the sample functions,
\[
 G(x):=\mathbb E_\xi[G(x;\xi)]\,,
 \qquad\text{where}\qquad
 G(x;\xi):=\nabla_y\,\mathbb P\bigl(C(x,y;\xi)=1\bigr)\big|_{y=x}\,.
\]
For $\varepsilon>0$, we say that $x$ is
\emph{$\varepsilon$-stationary with respect to comparisons} if
$\|G(x)\|\le\varepsilon$.
\end{definition}

To interpret $G(x)$, let $v$ be a unit vector. When the comparison
probability is differentiable at $y=x$, the oracle prefers $x+hv$ to $x$ with
probability $\tfrac12-h\langle G(x),v\rangle+o(h)$ as $h\to0$, uniformly in $v$.
Therefore, $\|G(x)\|$ is the largest rate at which a small move away from $x$
can increase the probability that the new point is preferred, and this rate is
attained in the direction of $-G(x)$. At an $\varepsilon$-stationary
point, every point $x+hv$ is preferred to $x$ with probability at most
$\tfrac12+\varepsilon h+o(h)$.
Under the oracles in \cref{def:exact-bounded-error}, by contrast,
$\mathbb P(C(x,x+hv)=1)$ jumps at $h=0$ whenever
$\langle\nabla f(x),v\rangle\ne0$, so the comparison probability has
no gradient at $y=x$.

Both $G(x)$ and $\varepsilon$-stationarity with respect to comparisons are
defined by the responses of the oracle alone. In black-box applications, such
as comparisons made by people or by language models, $f$ is not observed, and
these are the first-order notions that can be checked from comparisons. Under
\cref{def:smooth-oracle}, $G(x)$ is also a fixed positive multiple of
$\nabla f(x)$, as we show next.

\begin{lemma}
\label[lemma]{lem:comparison-gradient}
Suppose that the comparison oracle satisfies \cref{def:smooth-oracle}.
\begin{enumerate}[label=(\roman*),leftmargin=2.5em]
\item In deterministic optimization, if $f$ is differentiable at $x$, then
$G(x)$ exists and
\begin{equation}
 G(x)=\gamma\nabla f(x)\,.
 \label{eq:comparison-gradient-identity}
\end{equation}
\item In stochastic optimization, suppose that $f(\cdot;\xi)$ is
differentiable at $x$ for almost every $\xi$. Then $G(x)$ exists if and only
if $\mathbb E_\xi\|\nabla f(x;\xi)\|<\infty$, and in that case
$G(x)=\gamma\,\mathbb E_\xi\nabla f(x;\xi)$, which is $\gamma\nabla f(x)$ when
the gradient and the expectation commute.
\end{enumerate}
\end{lemma}
\begin{proof}
Part~(i) follows from~\eqref{eq:comparison-probability},
$\Phi'(0)=\gamma$, and the chain rule. For part~(ii), part~(i) applied to
$f(\cdot;\xi)$ gives $G(x;\xi)=\gamma\nabla f(x;\xi)$ for almost every $\xi$,
so $G(x)=\mathbb E_\xi[G(x;\xi)]$ exists if and only if
$\mathbb E_\xi\|\nabla f(x;\xi)\|<\infty$, and then
$G(x)=\gamma\,\mathbb E_\xi\nabla f(x;\xi)$.
\end{proof}

The factor $\gamma$, however, cannot be determined from comparisons.

\begin{proposition}
\label[proposition]{prop:gamma-nonidentifiability}
Let $\Phi$ satisfy \cref{def:smooth-oracle} with parameters
$(\gamma,M_\phi)$, and let $\tilde\Phi(t):=\Phi(\tfrac t2)$. For a
differentiable $f$, consider the oracle~\eqref{eq:comparison-probability}
with $\Phi$ for $f$ and with $\tilde\Phi$ for $2f$. The two problems have the
same comparison gradient $\gamma\nabla f(x)$ at $x$, but in both
problem~\eqref{eq:deterministic-problem} and problem~\eqref{eq:problem} with
$f(\cdot;\xi)=f$, any gradient estimator at $x$ has mean-squared error at least
$\tfrac14\|\nabla f(x)\|^2$ for one of them.
\end{proposition}
\begin{proof}
The function $\tilde\Phi$ satisfies \cref{def:smooth-oracle} with
parameters $(\tfrac\gamma2,\tfrac{M_\phi}{8})$, and $\tilde\Phi(2t)=\Phi(t)$.
Hence, by~\eqref{eq:comparison-probability}, every comparison returns $1$
with the same probability in the two problems, and,
by~\eqref{eq:comparison-gradient-identity}, their comparison gradients at $x$
are $\gamma\nabla f(x)$ and $\tfrac\gamma2\cdot2\nabla f(x)=\gamma\nabla f(x)$.
As in the proof of \cref{prop:deterministic-nonidentifiability}, $\hat g$ has the same
distribution in the two problems, whose gradients at $x$ are $\nabla f(x)$ and
$2\nabla f(x)$. As $\|g-a\|^2+\|g-2a\|^2\ge\tfrac12\|a\|^2$ for every
$g\in\mathbb R^d$, with $a=\nabla f(x)$, one of the two mean-squared errors is
at least $\tfrac14\|\nabla f(x)\|^2$.
\end{proof}

For the remainder of the paper, we consider smooth comparison oracles and
omit the term ``smooth'' when there is no ambiguity.

\FloatBarrier

\section{Estimating the comparison gradient}
\label{sec:estimation}
In this section, we estimate the comparison gradient $G(x)$ for
problems~\eqref{eq:deterministic-problem} and~\eqref{eq:problem}. Random-direction finite differences, which estimate the
gradient from function values at $x$ and at $x+hu$ for random directions $u$,
are standard in zeroth-order and derivative-free optimization
\citep{BerahasCaoChoromanskiScheinberg2022} and attain optimal rates in
stochastic convex optimization \citep{duchi2015optimal}. Their directions
may be uniform on the sphere \citep{FlaxmanKalaiMcMahan2005} or Gaussian
\citep{NesterovSpokoiny2017}. Since the comparison probability equals
$\tfrac12$ at $y=x$, the expected response of a comparison of $x$ and $x+hu$
is twice the difference of this probability between $x+hu$ and $x$. We
therefore extend these estimators to comparison oracles in a unified form,
the random-direction comparison estimator, which allows general spherically
symmetric directions (\cref{def:spherical-directions}). We bound its bias
and variance for problem~\eqref{eq:deterministic-problem} when the gradient
is Lipschitz and
when the Hessian is also Lipschitz (\cref{thm:det-gradient-estimator}), and,
as for function values
\citep{NesterovSpokoiny2017,AkhavanChzhenPontilTsybakov2024}, the bias
improves from order $h$ to order $h^2$ in the second case. We then extend
these bounds to problem~\eqref{eq:problem} (\cref{thm:gradient-estimator}).

We use a comparison oracle with parameters $(\gamma,M_\phi)$ from
\cref{def:smooth-oracle} and a differentiable function
$f:\mathbb R^d\to\mathbb R$, with the following smoothness conditions.
\begin{assumption}[Lipschitz gradient]\label{ass:function-smooth}
There exists $L_f\ge0$ such that $f$ is $L_f$-smooth:
\begin{equation}
\label{eq:function-smoothness}
        \left\lVert \nabla f(y)-\nabla f(x)\right\rVert\le L_f\left\lVert y-x\right\rVert\,,
        \quad\text{for all } x,y\in\mathbb{R}^d\,.
\end{equation}
\end{assumption}
For the sharper bound, we additionally assume a Lipschitz Hessian.
\begin{assumption}[Lipschitz Hessian]\label{ass:higher-order-smooth}
There exists $M_f\ge0$ such that $f$ is twice differentiable and has an
$M_f$-Lipschitz Hessian in the spectral norm $\|\cdot\|_{\mathrm{op}}$:
\begin{equation*}
        \left\lVert \nabla^2 f(y)-\nabla^2 f(x)\right\rVert_{\mathrm{op}}
        \le M_f\left\lVert y-x\right\rVert\,,
        \quad\text{for all } x,y\in\mathbb{R}^d\,.
\end{equation*}
\end{assumption}

\begin{estimatorbox}{Random-direction comparison estimator}
Let $x\in\mathbb R^d$ be an estimation point, $h>0$ a sampling radius, and
$1\le s\le d$ an integer. Draw $U=(u_1,\ldots,u_s)\sim\mathcal D_{d,s}$
of either type in \cref{def:spherical-directions} below.
For each direction $u_i$, draw a sample $\xi_i\sim\mathbb P_\xi$ and a
comparison variable $\zeta_i$, compare $x$ and $x+hu_i$ on the sample
function $f(\cdot;\xi_i)$, and set
\begin{equation}
 \hat g_h(x)=\frac{d}{2hs}
 \sum_{i=1}^s C_{\zeta_i}(x,x+hu_i;\xi_i)u_i\,.
 \label{eq:estimator}
\end{equation}
The samples and comparison variables are mutually independent and independent
of $U$. For problem~\eqref{eq:deterministic-problem}, $f(\cdot;\xi_i)=f$,
and we write
$C_{\zeta_i}(x,x+hu_i)$ for $C_{\zeta_i}(x,x+hu_i;\xi_i)$.
\end{estimatorbox}

\begin{definition}[Spherically symmetric directions]\label[definition]{def:spherical-directions}
Let $U=(u_1,\ldots,u_s)\sim\mathcal D_{d,s}$, where $\mathcal D_{d,s}$ is a
distribution on $(\mathbb R^d)^s$ and $1\le s\le d$.
We call $U$ a set of \emph{spherically symmetric directions} with a
dimension-independent parameter $K_{\mathcal D}\ge1$ if, for each
$i=1,\ldots,s$,
\begin{equation}
\label{eq:individual-spherical-symmetry}
 Qu_i\overset{\mathrm d}=u_i
 \quad\text{for any }Q\in\mathbb R^{d\times d}\text{ satisfying }Q^\top Q=I_d\,,
\end{equation}
\begin{equation}
\label{eq:individual-radius-conditions}
 \mathbb E\|u_i\|^2=1\,,
 \qquad\text{and}\qquad
 \mathbb E\exp\left(\frac{\|u_i\|^2}{K_{\mathcal D}^2}\right)\le2\,,
\end{equation}
and at least one of the following two conditions holds:
\begin{enumerate}[label=(\roman*),leftmargin=2.5em]
\item \emph{Independent directions}: $u_1,\ldots,u_s$ are mutually independent.
\item \emph{Orthonormal directions}: $u_1,\ldots,u_s$ are orthonormal almost
surely.
\end{enumerate}
\end{definition}
Condition~\eqref{eq:individual-spherical-symmetry} says that $Qu_i$ has the
same distribution as $u_i$ for every orthogonal matrix $Q$, and
condition~\eqref{eq:individual-radius-conditions} normalizes the directions
and bounds the tails of their lengths uniformly in the dimension. Many of the
random directions used in finite-difference methods and in random subspace
methods satisfy \cref{def:spherical-directions}, including the following
three.

\begin{example}[Gaussian directions]\label[example]{ex:gaussian-directions}
Let $u_1,\ldots,u_s$ be independent with distribution $N(0,\tfrac1dI_d)$.
Each $u_i$ is rotationally invariant, so
\eqref{eq:individual-spherical-symmetry} holds.
Also, $\mathbb E\|u_i\|^2=1$ and
$\mathbb E\exp(\tfrac14\|u_i\|^2)={(1-\tfrac1{2d})^{-d/2}}$ is at most $\sqrt2$, so
$U$ is of independent type in \cref{def:spherical-directions} with
$K_{\mathcal D}=2$.
\end{example}

\begin{example}[Uniform spherical directions]\label[example]{ex:spherical-directions}
Let $u_1,\ldots,u_s$ be independent and uniform on $\mathbb S^{d-1}$.
The uniform distribution is rotationally invariant, which implies
\eqref{eq:individual-spherical-symmetry}.
Since $\|u_i\|=1$, we have $\mathbb E\|u_i\|^2=1$ and
$\mathbb E\exp(\tfrac14\|u_i\|^2)=e^{1/4}<2$.
Therefore, $U$ is of independent type in \cref{def:spherical-directions}
with $K_{\mathcal D}=2$.
\end{example}

\begin{example}[Haar orthonormal directions]\label[example]{ex:haar-directions}
Let $U=(u_1,\ldots,u_s)$ be a Haar-distributed orthonormal $s$-frame.
Each $u_i$ is uniform on $\mathbb S^{d-1}$. Hence, as in
\cref{ex:spherical-directions}, $U$ is of orthonormal type in
\cref{def:spherical-directions} with $K_{\mathcal D}=2$.
\end{example}

Note that, for $d\ge2$, randomizing only coordinate directions and their
signs does not satisfy~\eqref{eq:individual-spherical-symmetry}. Similar
estimators with a single direction are used in the methods of
\citet{yue2009interactively} and \citet{kumagai2017regret}.

\subsection{Variance analysis}\label{sec:variance-analysis}
\begin{lemma}\label[lemma]{lem:variance-estimator}
Let $\hat g_h(x)$ be the estimator~\eqref{eq:estimator} with $U$ of either
type in \cref{def:spherical-directions}. Then
\begin{equation*}
\mathbb{E}\left[\left\lVert \hat g_h(x)-\mathbb{E}[\hat g_h(x)]\right\rVert^2\right]
\le \frac{d^2}{4h^2s}\,.
\end{equation*}
\end{lemma}
\begin{proof}
First suppose that $U$ is of orthonormal type. Using
$\mathbb E\|X-\mathbb EX\|^2\le\mathbb E\|X\|^2$, the orthonormality of
$u_1,\ldots,u_s$, $C_{\zeta_i}(x,x+hu_i)^2=1$, and $\|u_i\|=1$, we obtain
\begin{align*}
\mathbb{E}\left\lVert \hat g_h(x)-\mathbb{E}[\hat g_h(x)]\right\rVert^2
&\le \mathbb{E}\left\lVert \hat g_h(x)\right\rVert^2
\by{\eqref{eq:estimator}}{=}\frac{d^2}{4h^2s^2}\mathbb{E}\left[\sum_{i,j=1}^s
C_{\zeta_i}(x,x+hu_i)C_{\zeta_j}(x,x+hu_j)u_i^\top u_j\right]\\*
&=\frac{d^2}{4h^2s^2}\sum_{i=1}^s
\mathbb{E}\left[C_{\zeta_i}(x,x+hu_i)^2\left\lVert u_i\right\rVert^2\right]
=\frac{d^2}{4h^2s}\,.
\end{align*}
Now suppose that $U$ is of independent type. The vectors
$C_{\zeta_i}(x,x+hu_i)u_i$ are independent across $i$, so their
variances add. Bounding each variance by the second moment and using
$C_{\zeta_i}(x,x+hu_i)^2=1$, we obtain
\begin{align*}
\mathbb{E}\left\lVert \hat g_h(x)-\mathbb{E}[\hat g_h(x)]\right\rVert^2
&=\frac{d^2}{4h^2s^2}\sum_{i=1}^s
\mathbb{E}\left\lVert C_{\zeta_i}(x,x+hu_i)u_i
-\mathbb{E}[C_{\zeta_i}(x,x+hu_i)u_i]\right\rVert^2\\*
&\le\frac{d^2}{4h^2s^2}\sum_{i=1}^s
\mathbb{E}\left\lVert C_{\zeta_i}(x,x+hu_i)u_i\right\rVert^2
\by{\eqref{eq:individual-radius-conditions}}{=}\frac{d^2}{4h^2s}\,.\qedhere
\end{align*}
\end{proof}
The same calculation applies to problem~\eqref{eq:problem}, with $\zeta_i$
replaced by $(\xi_i,\zeta_i)$, since these pairs are independent across $i$
and independent of $U$.

\subsection{Bias analysis}\label{sec:bias-analysis}
To analyze the bias for problem~\eqref{eq:deterministic-problem}, fix
$x\in\mathbb R^d$ and define
\begin{equation}
 \phi_x(y):=\Phi(f(x+y)-f(x))\,.
 \label{eq:phi-x}
\end{equation}
By~\eqref{eq:comparison-probability},
$\mathbb E_\zeta[\tfrac{C_\zeta(x,x+y)}{2}\mid x,y]=\phi_x(y)$, which means that
$\phi_x$ is the centered probability that the oracle reports $f(x+y)>f(x)$.
Moreover, $\phi_x(0)=0$ and, by \cref{def:comparison-gradient},
$\nabla\phi_x(0)=G(x)$. Taking expectation in~\eqref{eq:estimator}, we
obtain
\begin{equation}
 \mathbb E[\hat g_h(x)]-G(x)
 =\frac1s\sum_{i=1}^s\left[
 \frac{d}{h}\mathbb E_U[(\phi_x(hu_i)-\phi_x(0))u_i]-G(x)\right]\,.
 \label{eq:deterministic-bias-averaging}
\end{equation}
Therefore, it suffices to bound the per-direction errors in brackets.
For the rest of this section, let $R_i:=\|u_i\|$ for each $i=1,\ldots,s$,
and let $w$ be uniform on $\mathbb B^d$ and independent of the full
batch $U$. The next lemma expresses each directional finite difference as
a weighted average of gradients over the unit ball.

\begin{lemma}\label[lemma]{lem:random-radius-ball-identity}
Under \cref{ass:function-smooth}, the function $\phi_x$ in
\eqref{eq:phi-x} satisfies, for each $i=1,\ldots,s$,
\begin{equation}
\label{eq:random-radius-ball-identity}
        \frac d h\mathbb{E}_{U\sim\mathcal D_{d,s}}[(\phi_x(hu_i)-\phi_x(0))u_i]
        =\mathbb{E}_{R_i,w}[R_i^2\nabla \phi_x(hR_iw)]\,.
\end{equation}
\end{lemma}

Before proving it, we recall the following identity.

\begin{lemma}[Lemma~2.1 of \citet{FlaxmanKalaiMcMahan2005}]
\label[lemma]{lem:ball-sphere-identity}
Let $g:\mathbb R^d\to\mathbb R$ be continuous and $\delta>0$, and let $w$ and
$q$ be uniform on $\mathbb B^d$ and $\mathbb S^{d-1}$, respectively. Then
$\hat g(y):=\mathbb E_w[g(y+\delta w)]$ is differentiable, and
$\nabla\hat g(y)=\tfrac d\delta\mathbb E_q[g(y+\delta q)q]$ for every
$y\in\mathbb R^d$.
\end{lemma}

\begin{proof}[Proof of \cref{lem:random-radius-ball-identity}]
Because $\Phi'$ is bounded by~\eqref{eq:Phi-bounds},
\cref{ass:function-smooth} implies that $\nabla\phi_x$ grows at most
linearly, so the expectations below are finite
by~\eqref{eq:individual-radius-conditions}. Fix $r>0$ and apply
\cref{lem:ball-sphere-identity} with $g=\phi_x$, $\delta=hr$, and $y=0$.
Since $\nabla\phi_x$ is continuous, the gradient of $\hat g$ at $0$ is
$\mathbb E_w[\nabla\phi_x(hrw)]$, and since $\mathbb E_qq=0$, we obtain
\[
 \mathbb E_w[\nabla\phi_x(hrw)]
 =\frac{d}{hr}\,\mathbb E_q[(\phi_x(hrq)-\phi_x(0))q]\,.
\]
By~\eqref{eq:individual-spherical-symmetry}, conditional on $R_i$, the
vector $u_i$ is uniform on the sphere of radius $R_i$, so it has the
distribution of $R_iq$. Setting $r=R_i$ and multiplying both sides by
$R_i^2$, we obtain, for $R_i>0$,
\[
 R_i^2\,\mathbb E_w[\nabla\phi_x(hR_iw)]
 =\frac dh\,\mathbb E\bigl[(\phi_x(hu_i)-\phi_x(0))u_i\bigm|R_i\bigr]\,.
\]
Both sides are zero when $R_i=0$. Taking expectations over $R_i$
proves~\eqref{eq:random-radius-ball-identity}.
\end{proof}

The next lemma bounds $\nabla\phi_x(y)-\nabla\phi_x(0)$.

\begin{lemma}
\label[lemma]{lem:local-gradient-bounds}
Under \cref{ass:function-smooth}, for any $y\in\mathbb R^d$,
\begin{equation}
 \|\nabla\phi_x(y)-\nabla\phi_x(0)\|
 \le\gamma L_f\|y\|
 +\frac52M_\phi\|\nabla f(x)\|^3\|y\|^2
 +M_\phi L_f^3\|y\|^5\,.
 \label{eq:local-gradient-bound-1}
\end{equation}
If \cref{ass:higher-order-smooth} also holds, then
\begin{equation}
 \|\nabla\phi_x(y)-\nabla\phi_x(0)-\gamma\nabla^2 f(x)y\|
 \le\frac{\gamma M_f}{2}\|y\|^2
 +\frac52M_\phi\|\nabla f(x)\|^3\|y\|^2
 +M_\phi L_f^3\|y\|^5\,.
 \label{eq:local-gradient-bound-2}
\end{equation}
\end{lemma}
\begin{proof}
By the chain rule,
\[
 \nabla\phi_x(y)=\Phi'(f(x+y)-f(x))\nabla f(x+y)\,,
\]
and in particular $\nabla\phi_x(0)=\gamma\nabla f(x)$, since $\Phi'(0)=\gamma$.
By \cref{ass:function-smooth},
\begin{equation}
 |f(x+y)-f(x)|\le\|y\|\left(\|\nabla f(x)\|+\frac{L_f\|y\|}{2}\right)\,,
 \qquad\text{and}\qquad \|\nabla f(x+y)\|\le\|\nabla f(x)\|+L_f\|y\|\,.
 \label{eq:basic-smoothness-bounds}
\end{equation}
Hence, by~\eqref{eq:Phi-bounds}
and~\eqref{eq:basic-smoothness-bounds},
\begin{equation}
\begin{aligned}
 &\|\nabla\phi_x(y)-\gamma \nabla f(x+y)\|
 =|\Phi'(f(x+y)-f(x))-\gamma|\,\|\nabla f(x+y)\|\\
 &\quad\le\frac{M_\phi}{2}\|y\|^2
 \bigl(\|\nabla f(x)\|+\tfrac{L_f\|y\|}{2}\bigr)^2
 (\|\nabla f(x)\|+L_f\|y\|)
 \le\frac52M_\phi\|\nabla f(x)\|^3\|y\|^2
 +M_\phi L_f^3\|y\|^5\,.
\end{aligned}
\label{eq:phi-x-gradient-error}
\end{equation}
For the last inequality, write $a=\|\nabla f(x)\|$ and $b=L_f\|y\|$, expand
the product, and use $a^2b\le\tfrac{2a^3+b^3}{3}$ and
$ab^2\le\tfrac{a^3+2b^3}{3}$. Since
\[
 \nabla\phi_x(y)-\nabla\phi_x(0)
 =[\nabla\phi_x(y)-\gamma\nabla f(x+y)]+\gamma[\nabla f(x+y)-\nabla f(x)]\,,
\]
the first term is bounded by~\eqref{eq:phi-x-gradient-error}
and the second by \cref{ass:function-smooth}, which proves
\eqref{eq:local-gradient-bound-1}.

Under \cref{ass:higher-order-smooth}, Taylor's theorem implies that
$\|\nabla f(x+y)-\nabla f(x)-\nabla^2 f(x)y\|
\le\tfrac{M_f}{2}\|y\|^2$. By the triangle inequality,
\begin{align*}
 \|\nabla\phi_x(y)-\nabla\phi_x(0)-\gamma\nabla^2 f(x)y\|
 &\le\|\nabla\phi_x(y)-\gamma\nabla f(x+y)\|
 +\gamma\|\nabla f(x+y)-\nabla f(x)-\nabla^2 f(x)y\|\\
 &\le\frac{\gamma M_f}{2}\|y\|^2
 +\frac52M_\phi\|\nabla f(x)\|^3\|y\|^2
 +M_\phi L_f^3\|y\|^5\,,
\end{align*}
where the last inequality uses~\eqref{eq:phi-x-gradient-error}
and the Taylor bound above. This proves
\eqref{eq:local-gradient-bound-2}.
\end{proof}

To state the bounds, we define the following constants:
\begin{equation}
 \kappa_1:=\gamma K_{\mathcal D}^3L_f\,,
 \quad \kappa_2:=\gamma K_{\mathcal D}^4M_f\,,
 \quad \kappa_3:=5M_\phi K_{\mathcal D}^4\,,
 \quad\text{and}\quad \kappa_4:=5M_\phi K_{\mathcal D}^7L_f^3\,.
\label{eq:kappa-definitions}
\end{equation}
Here $K_{\mathcal D}$ controls the tails of the direction lengths and equals
$2$ in \cref{ex:gaussian-directions,ex:spherical-directions,ex:haar-directions}.

\begin{lemma}
\label[lemma]{lem:per-direction-error}
Under \cref{ass:function-smooth}, for each $i=1,\ldots,s$,
\begin{equation}
 \left\|\frac{d}{h}\mathbb E_U[(\phi_x(hu_i)-\phi_x(0))u_i]-G(x)\right\|
 \le\kappa_1h+\kappa_3\|\nabla f(x)\|^3h^2+\kappa_4h^5\,.
 \label{eq:per-direction-error-1}
\end{equation}
If \cref{ass:higher-order-smooth} also holds, then
\begin{equation}
 \left\|\frac{d}{h}\mathbb E_U[(\phi_x(hu_i)-\phi_x(0))u_i]-G(x)\right\|
 \le\bigl(\kappa_2+\kappa_3\|\nabla f(x)\|^3\bigr)h^2+\kappa_4h^5\,.
 \label{eq:per-direction-error-2}
\end{equation}
\end{lemma}

\begin{proof}
For a fixed $i\in\{1,\ldots,s\}$, we first compute the radial moment bounds needed to apply
\cref{lem:local-gradient-bounds} in the weighted average
of \cref{lem:random-radius-ball-identity}.
For $q\in\{3,4,7\}$,
\[
 \mathbb E R_i^q
 \le K_{\mathcal D}^q\left(\frac{q}{2e}\right)^{q/2}
       \mathbb E\exp\Bigl(\frac{R_i^2}{K_{\mathcal D}^2}\Bigr)
 \by{\eqref{eq:individual-radius-conditions}}{\le}2\left(\frac{q}{2e}\right)^{q/2}K_{\mathcal D}^q\,.
\]
The first inequality uses
$t^{q/2}\le(\tfrac{q}{2e})^{q/2}e^t$ for $t\ge0$.
Since $\|w\|\le1$, it follows that
\begin{equation}
 \mathbb E[R_i^3\|w\|]\le K_{\mathcal D}^3\,,
 \quad \mathbb E[R_i^4\|w\|^2]\le2K_{\mathcal D}^4\,,
 \quad\text{and}\quad \mathbb E[R_i^7\|w\|^5]\le5K_{\mathcal D}^7\,.
 \label{eq:radial-moment-specializations}
\end{equation}
Since $\mathbb E R_i^2=1$ and $\nabla\phi_x(0)=G(x)$,
taking norms and using the triangle inequality, we obtain
\begin{align*}
 &\left\|\frac{d}{h}\mathbb E_U[(\phi_x(hu_i)-\phi_x(0))u_i]-G(x)\right\|
 \by{\cref{lem:random-radius-ball-identity}}{=}\left\|\mathbb E_{R_i,w}
 [R_i^2(\nabla\phi_x(hR_iw)-\nabla\phi_x(0))]\right\|\\
 &\quad\by{\eqref{eq:local-gradient-bound-1}}{\le}\gamma L_fh\mathbb E[R_i^3\|w\|]
 +\frac{5M_\phi}{2}\|\nabla f(x)\|^3h^2\mathbb E[R_i^4\|w\|^2]
 +M_\phi L_f^3h^5\mathbb E[R_i^7\|w\|^5]\,.
\end{align*}
Substituting~\eqref{eq:radial-moment-specializations}
and~\eqref{eq:kappa-definitions} proves~\eqref{eq:per-direction-error-1}.

Under \cref{ass:higher-order-smooth}, independence of $R_i,w$ and
$\mathbb E w=0$ imply
$\mathbb E[R_i^3\nabla^2 f(x)w]=0$.
Hence, by \cref{lem:random-radius-ball-identity} and $\mathbb E R_i^2=1$, we have
\begin{align*}
&\Bigl\|\frac{d}{h}\mathbb E_U[(\phi_x(hu_i)-\phi_x(0))u_i]-G(x)\Bigr\|
=\bigl\|\mathbb E_{R_i,w}[R_i^2(\nabla\phi_x(hR_iw)-\nabla\phi_x(0)-\gamma hR_i\nabla^2 f(x)w)]\bigr\|\\
&\quad\le\frac{\gamma M_fh^2}{2}\mathbb E[R_i^4\|w\|^2]
 +\frac{5M_\phi}{2}\|\nabla f(x)\|^3h^2\mathbb E[R_i^4\|w\|^2]
 +M_\phi L_f^3h^5\mathbb E[R_i^7\|w\|^5]\,,
\end{align*}
where the last inequality applies~\eqref{eq:local-gradient-bound-2}
with $y=hR_iw$ and the triangle inequality. Substituting
\eqref{eq:radial-moment-specializations} and~\eqref{eq:kappa-definitions}
proves~\eqref{eq:per-direction-error-2}.
\end{proof}

Combining \cref{lem:variance-estimator,lem:per-direction-error}
gives the following bounds for problem~\eqref{eq:deterministic-problem}.

\begin{theorem}\label{thm:det-gradient-estimator}
Let $\hat g_h(x)$ be the estimator~\eqref{eq:estimator} with $U$ of either
type in \cref{def:spherical-directions}, and let $\kappa_1,\ldots,\kappa_4$ be defined
in~\eqref{eq:kappa-definitions}.
Under \cref{ass:function-smooth},
\begin{equation*}
\begin{aligned}
 \|\mathbb E\hat g_h(x)-G(x)\|
 &\le\kappa_1h+\kappa_3\|\nabla f(x)\|^3h^2+\kappa_4h^5\,.
\end{aligned}
\end{equation*}
If \cref{ass:higher-order-smooth} also holds, then
\begin{equation*}
\begin{aligned}
 \|\mathbb E\hat g_h(x)-G(x)\|
 &\le\bigl(\kappa_2+\kappa_3\|\nabla f(x)\|^3\bigr)h^2+\kappa_4h^5\,.
\end{aligned}
\end{equation*}

The variance satisfies
\[
 \mathbb E\|\hat g_h(x)-\mathbb E\hat g_h(x)\|^2
 \le\frac{d^2}{4h^2s}\,.
\]
\end{theorem}
\begin{proof}
The variance bound is \cref{lem:variance-estimator}.
Taking norms in~\eqref{eq:deterministic-bias-averaging}, using the triangle
inequality, and applying \cref{lem:per-direction-error} to each
summand, we obtain the two bias bounds.
\end{proof}

\subsection{Extension to problem (S)}
\label{sec:stochastic}
We now extend \cref{thm:det-gradient-estimator} to
problem~\eqref{eq:problem}, where each comparison is made on a single sample
function.

\begin{theorem}\label{thm:gradient-estimator}
Let $\hat g_h(x)$ be the estimator~\eqref{eq:estimator} with $U$ of either
type in \cref{def:spherical-directions}, and let $\kappa_1,\ldots,\kappa_4$ be defined
in~\eqref{eq:kappa-definitions}. Suppose that $G(x)$ exists and that
$f(\cdot;\xi)$ satisfies \cref{ass:function-smooth} for almost every $\xi$.
Then
\begin{equation}
\begin{aligned}
 \|\mathbb E\hat g_h(x)-G(x)\|
 &\le\kappa_1h+\kappa_3\mathbb E_\xi\|\nabla f(x;\xi)\|^3h^2
 +\kappa_4h^5\,.
\end{aligned}
\label{eq:stochastic-estimator-bias}
\end{equation}
If, in addition, $f(\cdot;\xi)$ satisfies \cref{ass:higher-order-smooth} for
almost every $\xi$, then
\begin{equation}
\begin{aligned}
 \|\mathbb E\hat g_h(x)-G(x)\|
 &\le\bigl(\kappa_2+\kappa_3\mathbb E_\xi\|\nabla f(x;\xi)\|^3\bigr)h^2
 +\kappa_4h^5\,.
\end{aligned}
\label{eq:stochastic-estimator-bias-higher}
\end{equation}

The variance satisfies
\[
 \mathbb E\|\hat g_h(x)-\mathbb E\hat g_h(x)\|^2
 \le\frac{d^2}{4h^2s}\,.
\]
\end{theorem}
\begin{proof}
The variance bound is \cref{lem:variance-estimator}.

For the bias bounds, since averaging over $\xi$ does not preserve the form
$\tfrac12+\Phi(f(y)-f(x))$, we apply \cref{lem:per-direction-error}
to each sample function and then average. As
in~\eqref{eq:phi-x}, define
$\phi_x(y;\xi):=\Phi(f(x+y;\xi)-f(x;\xi))$.
By~\eqref{eq:comparison-probability},
$\mathbb E[C_{\zeta_i}(x,x+hu_i;\xi_i)\mid U,\xi_i]=2\phi_x(hu_i;\xi_i)$.
Since $|\phi_x|\le\tfrac12$, $\mathbb E\|u_i\|\le1$, and $\xi_i$ is
independent of $U$, Fubini's theorem gives
\[
 \mathbb E\hat g_h(x)
 =\frac{d}{hs}\sum_{i=1}^s\mathbb E_\xi\mathbb E_U[\phi_x(hu_i;\xi)u_i]\,.
\]
Since $G(x)$ exists, \cref{lem:comparison-gradient}(ii) gives
$G(x)=\gamma\mathbb E_\xi\nabla f(x;\xi)$. Subtracting it and using
$\phi_x(0;\xi)=0$ yield
\[
 \mathbb E\hat g_h(x)-G(x)
 =\frac1s\sum_{i=1}^s\mathbb E_\xi\!\left[
 \frac{d}{h}\mathbb E_U[(\phi_x(hu_i;\xi)-\phi_x(0;\xi))u_i]
 -\gamma\nabla f(x;\xi)\right]\,.
\]
For almost every $\xi$, $\gamma\nabla f(x;\xi)$ is the comparison gradient
of $f(\cdot;\xi)$, which satisfies \cref{ass:function-smooth}, and $U$ is
independent of $\xi$, so \cref{lem:per-direction-error}
applies to each vector in brackets. Taking norms and using the triangle
inequality give~\eqref{eq:stochastic-estimator-bias}, and the second bound of
\cref{lem:per-direction-error}
gives~\eqref{eq:stochastic-estimator-bias-higher} in the same way.
\end{proof}
The bias bounds of \cref{thm:gradient-estimator} depend on the samples only
through the third moment $\mathbb E_\xi\|\nabla f(x;\xi)\|^3$. By Jensen's
inequality,
\[
 \mathbb E_\xi\|\nabla f(x;\xi)\|^3
 \ge\bigl(\|\mathbb E_\xi\nabla f(x;\xi)\|^2
 +\mathbb E_\xi\|\nabla f(x;\xi)-\mathbb E_\xi\nabla f(x;\xi)\|^2\bigr)^{3/2}\,,
\]
so it reflects both the size of the mean sample gradient and the spread of the
sample gradients.

We now turn these bias and variance bounds into a number of comparisons.

\begin{definition}[Mean-squared error]\label[definition]{def:gradient-estimation-mse}
For a problem $f$ and an estimator $\hat g$ of $G(x)$ that uses $T$
comparisons, let
\[
 \operatorname{MSE}_T(\hat g;f):=\mathbb E\|\hat g-G(x)\|^2\,,
\]
where the expectation is over the randomness of the estimator, the
comparisons, and, for problem~\eqref{eq:problem}, the samples.
\end{definition}

Since the mean-squared error of an average of $T$ independent realizations is
the squared bias plus $\tfrac1T$ times the variance,
\cref{thm:det-gradient-estimator,thm:gradient-estimator} give the following.

\begin{corollary}
\label[corollary]{cor:gradient-comparison-complexity}
In the setting of \cref{thm:det-gradient-estimator} or
\cref{thm:gradient-estimator}, let $\hat g$ be the average of $T$
independent realizations of~\eqref{eq:estimator} with $s=1$, and write
$m_3:=\mathbb E_\xi\|\nabla f(x;\xi)\|^3$. For any given $\varepsilon>0$,
let $h>0$ be sufficiently small, satisfying
$\kappa_1h+\kappa_3m_3h^2+\kappa_4h^5\le\tfrac\varepsilon2$, or
$(\kappa_2+\kappa_3m_3)h^2+\kappa_4h^5\le\tfrac\varepsilon2$ when
\cref{ass:higher-order-smooth} also holds. Then for any integer
$T\ge\tfrac{d^2}{2h^2\varepsilon^2}$, it holds that
$\operatorname{MSE}_T(\hat g;f)\le\varepsilon^2$.
\end{corollary}

\begin{remark}
\label[remark]{rem:comparison-counts}
Once $h$ satisfies the condition, the bound on $T$ depends only on $d$, $h$,
and $\varepsilon$, so $T$ can be chosen without knowledge of the problem. For
small $\varepsilon$, the condition holds for all
$h\le\tfrac{\varepsilon}{4\kappa_1}$, which is of order $\varepsilon$, and the
largest of these yields $T_1:=\lceil\tfrac{8d^2\kappa_1^2}{\varepsilon^4}\rceil$. Under
\cref{ass:higher-order-smooth}, if $\kappa_2+\kappa_3m_3>0$, the second
condition holds for all $h\le\sqrt{\tfrac{\varepsilon}{4(\kappa_2+\kappa_3m_3)}}$,
which is of order $\sqrt\varepsilon$, and the largest of these yields
$T_2:=\lceil\tfrac{2d^2(\kappa_2+\kappa_3m_3)}{\varepsilon^3}\rceil$. The orders of
$T_1$ and $T_2$ are the same as those for estimating the gradient from noisy
function values \citep{AkhavanChzhenPontilTsybakov2024}.
\end{remark}

The bias bounds of \cref{thm:det-gradient-estimator,thm:gradient-estimator}
involve the size of the gradient and the third moment
$\mathbb E_\xi\|\nabla f(x;\xi)\|^3$, respectively. In optimization with
stochastic gradients, neither is usually assumed to be bounded, since a bound
on the variance of the stochastic gradients suffices for estimating the
gradient. With comparisons, however, a bound on the size of the gradient is
necessary, even for problem~\eqref{eq:deterministic-problem}. For example,
consider, in one dimension, the problems $f_a(z):=az$ with $a\in\mathbb R$, of
the form~\eqref{eq:deterministic-problem}. They have $\nabla f_a(0)=a$ and
zero Hessian, so both smoothness assumptions hold with $L_f=M_f=0$. As
instances of problem~\eqref{eq:problem} with $f_a(z;\xi):=f_a(z)$, they also
have zero variance. Let $\mathbb P_a$ denote probability under $f_a$.

\begin{proposition}
\label[proposition]{prop:bounded-variance-insufficient}
Let $T\ge1$ be an integer, and let $\hat g\in\mathbb R$ be any estimator of
$G(0)=\gamma a$ that uses $T$ comparisons on the problems $f_a$, at any points
chosen possibly adaptively, and that may know $\Phi$ but not $a$. Then
$\sup_{a\in\mathbb R}\operatorname{MSE}_T(\hat g;f_a)=+\infty$.
\end{proposition}
\begin{proof}
Fix the random seed of $\hat g$, which does not depend on $a$. Then the
first query is fixed, each later query is determined by the previous
responses, and the output is determined by all $T$ responses. Hence
$\hat g$ takes at most $2^T$ values, wherever the queries lie.

Let $R>0$ and $a_j:=\tfrac{4Rj}{\gamma}$ for $j=1,\ldots,2^{T+1}$, so that
$\gamma a_j=4Rj$. Each of these values is within distance $R$ of at most one
$\gamma a_j$, so averaging over $j$ and over the seed gives
$2^{-T-1}\sum_{j=1}^{2^{T+1}}\mathbb P_{a_j}(|\hat g-\gamma a_j|<R)\le\tfrac12$.
Hence $\mathbb P_{a_j}(|\hat g-\gamma a_j|\ge R)\ge\tfrac12$, and so
$\operatorname{MSE}_T(\hat g;f_{a_j})\ge\tfrac{R^2}{2}$, for some $j$. Letting
$R$ grow completes the proof.
\end{proof}

\section{How many comparisons are necessary}
\label{sec:gradient-lower-bound}

This section asks how many comparisons are necessary for estimating
first-order information, namely the comparison gradient $G(x)$ at a given
point $x\in\mathbb R^d$. We call an estimator \emph{comparison-based} if it
makes $T$ comparisons as follows. It draws a random seed $\omega$ before the
first comparison. At round $t$, it chooses the pair $(x_t,y_t)$ as a
measurable function of $\omega$ and the previous responses
$b_1,\ldots,b_{t-1}$, and the oracle returns the response
$b_t:=C_{\zeta_t}(x_t,y_t)\in\{-1,1\}$. After $T$ rounds, it outputs $\hat g$,
a measurable function of $\omega$ and $b_1,\ldots,b_T$.

For $L_f,M_f\ge0$, let $\mathcal F_1(L_f)$ denote the class of functions
$f:\mathbb R^d\to\mathbb R$ that satisfy \cref{ass:function-smooth} with
constant $L_f$, and let $\mathcal F_2(L_f,M_f)$ denote the subclass of those
that also satisfy \cref{ass:higher-order-smooth} with constant $M_f$.

\begingroup
\begin{theorem}
\label{thm:comparison-complexity-lower-bounds}
There exist absolute constants $c_1,c_2,c_3,c_4>0$ such that, for any
$d\ge1$, any comparison oracle with parameters $(\gamma,M_\phi)$, and any
comparison-based estimator with output $\hat g$, the following statements
hold for problem~\eqref{eq:deterministic-problem}.
\begin{enumerate}[label=(\roman*),leftmargin=2.5em]
\item Let $L_f,B,\varepsilon>0$ satisfy $\varepsilon\le\tfrac B4$ and
$\varepsilon\le c_1(\tfrac{\gamma^5L_f^2}{M_\phi})^{1/4}$.
If $\operatorname{MSE}_T(\hat g;f)\le\varepsilon^2$ for every
$f\in\mathcal F_1(L_f)$ with $\|G(x)\|\le B$, then
\begin{equation}
T\ge
c_2\frac{\gamma^2d^2L_f^2}{\varepsilon^4}\,.
\label{eq:first-order-gradient-lower-bound}
\end{equation}

\item Let $L_f,M_f,B,\varepsilon>0$ satisfy $\varepsilon\le\tfrac B4$,
$\varepsilon\le c_3\tfrac{\gamma L_f^2}{M_f}$, and
$\varepsilon^3\le c_3^3\tfrac{\gamma^4M_f}{M_\phi}$.
If $\operatorname{MSE}_T(\hat g;f)\le\varepsilon^2$ for every
$f\in\mathcal F_2(L_f,M_f)$ with $\|G(x)\|\le B$, then
\begin{equation}
T\ge
c_4\frac{\gamma d^2M_f}{\varepsilon^3}\,.
\label{eq:second-order-gradient-lower-bound}
\end{equation}
\end{enumerate}
\end{theorem}
\endgroup
Problems of the form~\eqref{eq:deterministic-problem} are instances of
problem~\eqref{eq:problem} whose sample functions do not depend on the
sample, so the lower bounds also hold for problem~\eqref{eq:problem}. For
fixed $B$ and all sufficiently small $\varepsilon$, part~(i) matches $T_1$ in
\cref{rem:comparison-counts} up to a constant factor, and part~(ii) matches
$T_2$ in its dependence on $d$ and $\varepsilon$. Hence the random-direction
comparison estimator is optimal in $d$ and $\varepsilon$ under
\cref{ass:function-smooth} and under \cref{ass:higher-order-smooth}.

Concurrently with this work, \citet{LuEtAl2026DuelingFeedback} use a similar
estimator in online convex optimization, where $\Phi'(0)$ is known. Since
$\Phi'(0)$ cannot be determined from comparisons
(\cref{prop:gamma-nonidentifiability}), our estimator targets $G(x)$, and
\cref{cor:gradient-comparison-complexity,thm:comparison-complexity-lower-bounds}
give the number of comparisons needed to estimate it at a given point.

In a revision posted in October 2026,
\citet{ChaytiElBakkaliSaadi2026NoisyComparisons} added a lower bound of order
$d^2$ on the number of comparisons needed to find a Goldstein stationary point
when $\Phi$ is known, with a different hard family and a different proof. \Cref{thm:comparison-complexity-lower-bounds} instead
bounds the mean-squared error of estimating $G(x)$ at a given point, holds for
any $\Phi$ in the model, and matches the upper bounds under both smoothness
conditions.

\subsection{A hard family}
\label{subsec:hard-family}

For $R,\rho>0$, the \emph{hard family} consists of the functions $f_\theta$
with $\lVert\theta\rVert\le R$, where
\begin{equation}
        \Psi_\rho(z):=\frac{z}{1+\left\lVert z\right\rVert^2/\rho^2}\,,
        \qquad\text{and}\qquad
        f_\theta(z):=\left\langle \theta,\Psi_\rho(z-x)\right\rangle\,,
         \qquad\text{for all }z\in\mathbb{R}^d\,.
\label{eq:hard-family}
\end{equation}
\begin{wrapfigure}{r}{0.48\textwidth}
\centering
\vspace{-\baselineskip}
\includegraphics[width=\linewidth]{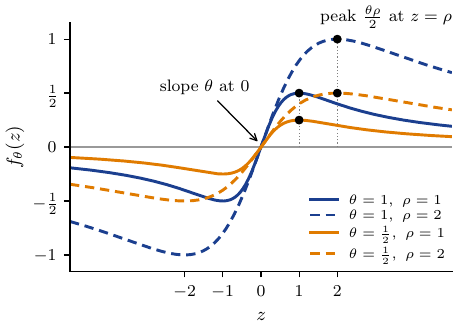}
\caption{The family~\eqref{eq:hard-family} with $d=1$ and $x=0$.}
\label{fig:hard-family}
\vspace{-0.3\baselineskip}
\end{wrapfigure}
Equivalently, $f_\theta$ is the linear function $\langle\theta,z-x\rangle$
multiplied by the factor $\tfrac{1}{1+\lVert z-x\rVert^2/\rho^2}$, as
\cref{fig:hard-family} shows in one dimension. Near $x$ this factor is close
to one, so $\nabla f_\theta(x)=\theta$, and the comparison gradient of
$f_\theta$ at $x$ is $\gamma\theta$ by~\eqref{eq:comparison-gradient-identity}.
Estimating $G(x)$ for this family is therefore estimating $\theta$. Far from
$x$, the factor brings $f_\theta$ back toward $f_\theta(x)=0$ for every
$\theta$, so only queried points within distance of order $\rho$ from $x$
provide information about $\theta$. Moreover, $f_\theta$ takes values
between $-\tfrac{R\rho}2$ and $\tfrac{R\rho}2$, so all comparison gaps are at
most $R\rho$, and when $R\rho$ is small, every comparison provides little
information about $\theta$, wherever the queried points lie.

For the fixed estimator, write $\mathbb P_\theta$ and $\mathbb E_\theta$
for the law and expectation of the record
$\mathcal T_T:=(\omega,(x_t,y_t,b_t)_{t=1}^T)$ when the objective is
$f_\theta$. Since $G(x)=\gamma\theta$ for $f_\theta$,
\cref{def:gradient-estimation-mse} yields
\begin{equation}
\sup_{\left\lVert \theta\right\rVert\le R}
\operatorname{MSE}_T(\hat g;f_\theta)
=\sup_{\left\lVert \theta\right\rVert\le R}
\mathbb{E}_\theta\left\lVert \hat g-\gamma\theta\right\rVert^2\,.
\label{eq:hard-family-worst-case-mse}
\end{equation}
By~\eqref{eq:hard-family-worst-case-mse}, it suffices to lower-bound the
worst-case mean-squared error for estimating $\gamma\theta$ from the
comparison responses.

We prove \cref{thm:comparison-complexity-lower-bounds} in two steps. First, we lower-bound the worst-case
mean-squared error on
this family in terms of the number of comparisons and the parameters $R$
and $\rho$. We then relate $R$ and $\rho$ to the smoothness bounds and choose
them so that the family belongs to each
problem class, obtaining lower bounds on the number of comparisons.

\subsection{A mean-squared error lower bound}
\label{subsec:mse-lower-bound}

For the hard family in~\eqref{eq:hard-family}, the function-value gap is
$f_\theta(y_t)-f_\theta(x_t)
=\langle\theta,\Psi_\rho(y_t-x)-\Psi_\rho(x_t-x)\rangle$.
Let $v_t:=\Psi_\rho(y_t-x)-\Psi_\rho(x_t-x)$. Since
$1+\tfrac{\lVert z\rVert^2}{\rho^2}\ge\tfrac{2\lVert z\rVert}{\rho}$, we have
$\lVert\Psi_\rho(z)\rVert\le\tfrac\rho2$, so $\lVert v_t\rVert\le\rho$. Hence,
\begin{equation}
        f_\theta(y_t)-f_\theta(x_t)
        =\left\langle \theta,v_t\right\rangle\,,
        \qquad\text{and}\qquad
        \left\lVert v_t\right\rVert\le\rho\,,
         \qquad\text{for }t=1,\ldots,T\,.
\label{eq:linear-comparison-gap}
\end{equation}

Throughout this subsection, let $R,\rho>0$ satisfy
$R\rho\le\tfrac{1}{8\gamma}$ and $\tfrac{M_\phi R^2\rho^2}{2}\le\gamma$. We
bound the Fisher information that the comparisons carry about $\theta$, and the
van Trees inequality turns this bound into a lower bound on the mean-squared
error.

For the fixed estimator, let $\mathcal H_{t-1}$ denote the $\sigma$-algebra
generated by $\omega$ and $b_1,\ldots,b_{t-1}$. Each query is determined by
the seed and the previous responses, so $v_t$ is fixed given
$\mathcal H_{t-1}$, by a rule that does not depend on $\theta$. Its
distribution may nevertheless depend on $\theta$ through the previous
responses. By \cref{def:smooth-oracle} and~\eqref{eq:linear-comparison-gap},
\begin{equation}
 \mathbb P_\theta(b_t\mid\mathcal H_{t-1})
 =\frac12+b_t\Phi(\langle\theta,v_t\rangle)\,.
\label{eq:conditional-comparison-probabilities}
\end{equation}
Given the seed, the probability of the observed responses is therefore
\[
 p_\theta(\mathcal T_T):=\prod_{t=1}^T
 \Bigl(\frac12+b_t\Phi(\langle\theta,v_t\rangle)\Bigr)\,.
\]
Since the distribution of $\omega$ does not depend on $\theta$, we use it as
the likelihood of the record $\mathcal T_T=(\omega,(x_t,y_t,b_t)_{t=1}^T)$ of
the seed, the queries, and the responses. It is continuously differentiable
in $\theta$, and its score has mean zero because the probabilities of the two
responses sum to one. The Fisher information matrix of $\mathcal T_T$ with
respect to $\theta$ is
\begin{equation}
I_T(\theta):=\mathbb E_\theta\!\left[
\nabla_\theta\log p_\theta(\mathcal T_T)
\bigl(\nabla_\theta\log p_\theta(\mathcal T_T)\bigr)^\top
\right]\,.
\label{eq:comparison-fisher-information}
\end{equation}
We use the following form of the van Trees inequality.
\begin{lemma}[van Trees inequality, inequality~(3) of \citet{GillLevit1995}]
\label[lemma]{lem:van-trees}
Let an observation $X$ have density $p_\theta$ with respect to a measure
$\nu$, for $\theta$ in a closed interval $\Theta\subset\mathbb R$, and let
$\lambda$ be a prior density on $\Theta$. Suppose that $\lambda$ and
$\theta\mapsto p_\theta(X)$, for $\nu$-almost all $X$, are absolutely
continuous, that $\lambda$ converges to zero at the endpoints of $\Theta$, and
that $\mathbb E_\theta[\partial_\theta\log p_\theta(X)]=0$. Then any
estimator $\hat\theta$ of $\theta$ based on $X$ satisfies
\[
 \mathbb E_{\theta\sim\lambda}\mathbb E_\theta(\hat\theta-\theta)^2
 \ge\frac{1}{\mathbb E_{\theta\sim\lambda}[I(\theta)]+J(\lambda)}\,,
\]
where $\theta$ is drawn from $\lambda$ and then $X$ from $p_\theta$. Here
$I(\theta):=\mathbb E_\theta[(\partial_\theta\log p_\theta(X))^2]$ is the
Fisher information for $\theta$, and
$J(\lambda):=\mathbb E_{\theta\sim\lambda}[(\partial_\theta\log\lambda(\theta))^2]$ is
that for a location parameter in $\lambda$.
\end{lemma}

For the hard family, the inequality takes the following form.

\begin{lemma}
\label[lemma]{lem:mse-fisher-information}
Any comparison-based estimator satisfies
\begin{equation}
 \sup_{\lVert\theta\rVert\le R}\operatorname{MSE}_T(\hat g;f_\theta)
 \ge\frac{\gamma^2d^2}{\sup_{\lVert\theta\rVert\le R}\operatorname{tr}I_T(\theta)
 +10d^2/R^2}\,.
\label{eq:mse-fisher-information}
\end{equation}
\end{lemma}
\begin{proof}
Let $a:=\tfrac{R}{\sqrt d}$, so that $\lVert\theta\rVert\le R$ for
$\theta\in[-a,a]^d$. Take the prior density
$\lambda(t):=\tfrac1a\cos^2(\tfrac{\pi t}{2a})$ on $[-a,a]$. It is
continuously differentiable and vanishes at $\pm a$, and
$\lambda'(t)=-\tfrac{\pi}{2a^2}\sin(\tfrac{\pi t}{a})$, so
\[
 J(\lambda)=\int_{-a}^a\frac{\lambda'(t)^2}{\lambda(t)}\,dt
 =\int_{-a}^a\frac{\pi^2}{a^3}\sin^2\Bigl(\frac{\pi t}{2a}\Bigr)\,dt
 =\frac{\pi^2}{a^2}\,.
\]
Let $\lambda_d(\theta):=\prod_{i=1}^d\lambda(\theta_i)$ on $[-a,a]^d$.

Fix $i$ and the other components $\theta_j$, $j\ne i$, in $[-a,a]$, and
regard $\theta_i\in[-a,a]$ as the parameter. By the properties of
$p_\theta$ stated before~\eqref{eq:comparison-fisher-information}, the
assumptions of \cref{lem:van-trees} hold, with $\nu$ the product of the
distribution of $\omega$ and counting measure on $\{-1,1\}^T$, and the Fisher
information about $\theta_i$ is $[I_T(\theta)]_{ii}$. Applying
\cref{lem:van-trees} to the estimator $\tfrac{\hat g_i}{\gamma}$ of $\theta_i$ gives
\[
 \int_{-a}^a\mathbb E_\theta\Bigl(\frac{\hat g_i}{\gamma}-\theta_i\Bigr)^2
 \lambda(\theta_i)\,d\theta_i
 \ge\frac{1}{\int_{-a}^a[I_T(\theta)]_{ii}\lambda(\theta_i)\,d\theta_i
 +\pi^2/a^2}\,.
\]
We average over the other components with density
$\prod_{j\ne i}\lambda(\theta_j)$, use the convexity of $s\mapsto\tfrac1s$, and then
sum over $i$, using the convexity again:
\begin{align*}
 \mathbb E_{\theta\sim\lambda_d}\mathbb E_\theta
 \bigl\|\tfrac{\hat g}{\gamma}-\theta\bigr\|^2
 &\ge\sum_{i=1}^d\frac{1}{\mathbb E_{\theta\sim\lambda_d}[I_T(\theta)]_{ii}
 +\pi^2/a^2}
 \ge\frac{d^2}{\mathbb E_{\theta\sim\lambda_d}[\operatorname{tr}I_T(\theta)]
 +d\pi^2/a^2}\\
 &\ge\frac{d^2}{\sup_{\lVert\theta\rVert\le R}\operatorname{tr}I_T(\theta)
 +10d^2/R^2}\,,
\end{align*}
where the last step uses $\tfrac{d\pi^2}{a^2}=\tfrac{\pi^2d^2}{R^2}\le\tfrac{10d^2}{R^2}$. Since
$\lambda_d$ is supported on $[-a,a]^d$, multiplying by $\gamma^2$
and~\eqref{eq:hard-family-worst-case-mse} give~\eqref{eq:mse-fisher-information}.
\end{proof}

It remains to bound the Fisher information of the record.

\begin{lemma}
\label[lemma]{lem:adaptive-information-bound}
For any comparison-based estimator and all $\lVert\theta\rVert\le R$, we have
$\operatorname{tr}I_T(\theta)\le22\gamma^2\rho^2T$.
\end{lemma}
\begin{proof}
We bound the information in one comparison and then accumulate these bounds
under adaptive querying. Cauchy--Schwarz, $\|\theta\|\le R$,
and~\eqref{eq:linear-comparison-gap} give
$|\langle\theta,v_t\rangle|\le\lVert\theta\rVert\lVert v_t\rVert\le R\rho$.
Let $S_t(\theta):=\nabla_\theta\log\mathbb P_\theta(b_t\mid\mathcal H_{t-1})$
be the conditional score in round $t$, where the derivative holds the
realized history and comparison response fixed.

For any $\eta$ between zero and $\langle\theta,v_t\rangle$,
\eqref{eq:Phi-bounds}, this bound, and
$\tfrac{M_\phi R^2\rho^2}{2}\le\gamma$ imply that
$|\Phi'(\eta)-\gamma|\le\tfrac{M_\phi}{2}\eta^2\le\tfrac{M_\phi}{2}R^2\rho^2\le\gamma$,
so $|\Phi'(\eta)|\le2\gamma$ on this segment. Integrating this bound
from zero to $\left\langle \theta,v_t\right\rangle$, with $\Phi(0)=0$ and
$|\langle\theta,v_t\rangle|\le R\rho\le\tfrac{1}{8\gamma}$, yields
\begin{equation}
\left|\Phi(\left\langle \theta,v_t\right\rangle)\right|\le\frac14\,,
\qquad\text{and}\qquad
\left|\Phi'(\left\langle \theta,v_t\right\rangle)\right|\le2\gamma\,.
\label{eq:one-comparison-local-bounds}
\end{equation}
By~\eqref{eq:conditional-comparison-probabilities}, the conditional score is
\[
 S_t(\theta)=\frac{b_t\Phi'(\langle\theta,v_t\rangle)}
 {\frac12+b_t\Phi(\langle\theta,v_t\rangle)}\,v_t\,.
\]
Multiplying the squared norms of the scores of the two responses by their
probabilities, we obtain
\begin{equation}
\begin{aligned}
 \mathbb E_\theta[\|S_t(\theta)\|^2\mid\mathcal H_{t-1}]
 &=\Phi'(\langle\theta,v_t\rangle)^2\|v_t\|^2
 \left(\frac{1}{\frac12-\Phi(\langle\theta,v_t\rangle)}
 +\frac{1}{\frac12+\Phi(\langle\theta,v_t\rangle)}\right)\\
 &=\frac{\Phi'(\langle\theta,v_t\rangle)^2}
 {\frac14-\Phi(\langle\theta,v_t\rangle)^2}\,\|v_t\|^2
 \le22\gamma^2\|v_t\|^2\,.
\end{aligned}
\label{eq:per-round-information-bound}
\end{equation}
The inequality uses~\eqref{eq:one-comparison-local-bounds}:
$\Phi'(\langle\theta,v_t\rangle)^2\le4\gamma^2$ and
$\tfrac14-\Phi(\langle\theta,v_t\rangle)^2\ge\tfrac{3}{16}$, which shows that the
numerical factor is at most $\tfrac{64}{3}<22$.

To accumulate these bounds, recall that $p_\theta(\mathcal T_T)$ is the
product of the conditional probabilities
in~\eqref{eq:conditional-comparison-probabilities}. Taking the gradient of its
logarithm, we obtain
$\nabla_\theta\log p_\theta(\mathcal T_T)=\sum_{t=1}^T S_t(\theta)$.

Direct summation over $b_t\in\{-1,1\}$ shows that
$\mathbb E_\theta[S_t(\theta)\mid\mathcal H_{t-1}]=0$. For $r<t$, the score
$S_r(\theta)$ is $\mathcal H_{t-1}$-measurable, so conditioning on
$\mathcal H_{t-1}$ gives
\[
 \mathbb E_\theta\langle S_t(\theta),S_r(\theta)\rangle
 =\mathbb E_\theta\bigl\langle\mathbb E_\theta[S_t(\theta)\mid\mathcal H_{t-1}],
 S_r(\theta)\bigr\rangle=0\,.
\]
Therefore, the cross terms vanish when we expand the
squared transcript score:
\[
 \operatorname{tr}I_T(\theta)
 \by{\eqref{eq:comparison-fisher-information}}{=}\mathbb E_\theta\Bigl\|\sum_{t=1}^T S_t(\theta)\Bigr\|^2
 =\sum_{t=1}^T\mathbb E_\theta\|S_t(\theta)\|^2
 \by{\eqref{eq:per-round-information-bound}}{\le}22\gamma^2\sum_{t=1}^T\mathbb E_\theta\|v_t\|^2
 \by{\eqref{eq:linear-comparison-gap}}{\le}22\gamma^2\rho^2T\,.\qedhere
\]
\end{proof}

Combining the two lemmas gives the mean-squared error lower bound used in
\cref{subsec:lower-bound-proofs}.

\begin{lemma}\label[lemma]{lem:worst-case-mse-lower-bound}
Any comparison-based estimator satisfies
\begin{equation}
        \sup_{\left\lVert \theta\right\rVert\le R}
        \operatorname{MSE}_T(\hat g;f_\theta)
        \ge
        \frac{\gamma^2d^2}
        {22\gamma^2\rho^2T+10d^2/R^2}\,.
\label{eq:worst-case-mse-lower-bound}
\end{equation}
\end{lemma}
\begin{proof}
Substituting the bound of \cref{lem:adaptive-information-bound}
into~\eqref{eq:mse-fisher-information} gives~\eqref{eq:worst-case-mse-lower-bound}.
\end{proof}

\subsection{Proof of the lower bounds}
\label{subsec:lower-bound-proofs}

\begin{lemma}\label[lemma]{lem:hard-family-smoothness}
For any $R,\rho>0$ and $\lVert\theta\rVert\le R$, the function $f_\theta$
satisfies \cref{ass:function-smooth} with constant $\tfrac{3R}{\rho}$ and
\cref{ass:higher-order-smooth} with constant $\tfrac{13R}{\rho^2}$.
\end{lemma}
We use this lemma to choose $R$ and $\rho$ below, and we prove it in
\hyperref[app:proofs]{the appendix}.

\begin{proof}[Proof of \cref{thm:comparison-complexity-lower-bounds}]
Let $\mathcal F$ denote $\mathcal F_1(L_f)$ in part~(i) and
$\mathcal F_2(L_f,M_f)$ in part~(ii). Let $\hat g$ be a comparison-based
estimator such that $\operatorname{MSE}_T(\hat g;f)\le\varepsilon^2$ for every
$f\in\mathcal F$ with $\|G(x)\|\le B$. We apply
\cref{lem:worst-case-mse-lower-bound} to the hard family with
$R:=\tfrac{4\varepsilon}{\gamma}$. Every $f_\theta$ in this family satisfies
$\|G(x)\|=\gamma\|\theta\|\le4\varepsilon\le B$. The conditions
$R\rho\le\tfrac1{8\gamma}$ and $\tfrac{M_\phi R^2\rho^2}{2}\le\gamma$ of
\cref{subsec:mse-lower-bound} hold if and only if
$R\rho\le\min\{\tfrac1{8\gamma},\sqrt{\tfrac{2\gamma}{M_\phi}}\}$. Suppose that this
inequality holds and that the hard family is contained in $\mathcal F$. Then
\cref{lem:worst-case-mse-lower-bound} yields
\begin{equation*}
\varepsilon^2
\ge\frac{\gamma^2d^2}{22\gamma^2\rho^2T+10d^2/R^2}
=\frac{\gamma^2d^2}{22\gamma^2\rho^2T+5\gamma^2d^2/(8\varepsilon^2)}\,,
\end{equation*}
and solving this inequality for $T$ gives
$T\ge\tfrac{3d^2}{176\rho^2\varepsilon^2}$.
It remains to choose $\rho$ in each part so that the hard family lies in
$\mathcal F$ and $R\rho$ satisfies the inequality above.

For part~(i), we choose $\rho:=\tfrac{3R}{L_f}$. By
\cref{lem:hard-family-smoothness}, every $f_\theta$ then satisfies
\cref{ass:function-smooth} with constant $L_f$, so the hard family lies in
$\mathcal F_1(L_f)$. Moreover, $R\rho=\tfrac{48\varepsilon^2}{\gamma^2L_f}$,
so $R\rho$ satisfies the inequality above when $\varepsilon$ is at most a
small absolute constant times
$(\tfrac{\gamma^5L_f^2}{M_\phi})^{1/4}$, because
$M_\phi\ge\tfrac{27}{8}\gamma^3$ gives
$(\tfrac{\gamma^5L_f^2}{M_\phi})^{1/4}\le\sqrt{\gamma L_f}$. This is
the assumption of part~(i). Substituting
$\rho=\tfrac{12\varepsilon}{\gamma L_f}$
into this bound on $T$ gives
\eqref{eq:first-order-gradient-lower-bound}.

For part~(ii), we choose $\rho:=\sqrt{\tfrac{13R}{M_f}}$. By
\cref{lem:hard-family-smoothness}, every $f_\theta$ then satisfies
\cref{ass:higher-order-smooth} with constant $M_f$, and it also satisfies
\cref{ass:function-smooth} with constant $L_f$ if $\tfrac{3R}{\rho}\le L_f$,
that is, if $\varepsilon\le\tfrac{13\gamma L_f^2}{36M_f}$. Moreover,
$(R\rho)^2=\tfrac{832\varepsilon^3}{\gamma^3M_f}$, so $R\rho$ satisfies the
inequality above when $\varepsilon^3$ is at most a small absolute constant
times $\tfrac{\gamma^4M_f}{M_\phi}$, because $M_\phi\ge\tfrac{27}{8}\gamma^3$
gives $\tfrac{\gamma^4M_f}{M_\phi}\le\gamma M_f$. These are the
assumptions of part~(ii). Substituting
$\rho^2=\tfrac{52\varepsilon}{\gamma M_f}$
into this bound on $T$ gives
\eqref{eq:second-order-gradient-lower-bound}.
\end{proof}

\section{Comparison-based optimization}
\label{sec:optimization-application}

Recall from \cref{def:comparison-gradient} that $x$ is
$\varepsilon$-stationary with respect to comparisons if
$\|G(x)\|\le\varepsilon$. We now use the random-direction comparison
estimator in a stochastic gradient method to find such a point.

We first recall a general framework for stochastic gradient methods with
inexact gradients. Let $\psi:\mathbb R^d\to\mathbb R$ be $L_\psi$-smooth with
$\psi_*:=\inf_x\psi(x)>-\infty$. Consider the stochastic gradient method
$x_{k+1}=x_k-\alpha\tilde g_k$, where $\tilde g_k$ estimates
$\nabla\psi(x_k)$ with new randomness at each iteration and, given the
history $\mathcal F_k$,
\begin{equation}
 \|\mathbb E[\tilde g_k\mid\mathcal F_k]-\nabla\psi(x_k)\|^2\le\epsilon_g^2\,,
 \qquad\text{and}\qquad
 \mathbb E[\|\tilde g_k-\mathbb E[\tilde g_k\mid\mathcal F_k]\|^2\mid\mathcal F_k]\le\sigma_c^2\,.
\label{eq:oracle-conditions}
\end{equation}
\citet[Theorem~4.18(1)]{berahas_cao_scheinberg_oracles_review}, with
$\kappa_g=\sigma_v=0$ and $\tau=\tfrac{\sqrt3}2$, shows that if
$\alpha\le\tfrac1{L_\psi}$ and
$\epsilon_g^2+\alpha L_\psi\sigma_c^2\le\tfrac{3\varepsilon^2}8$, then the
first iteration at which $\|\nabla\psi(x_k)\|\le\varepsilon$ satisfies
\begin{equation}
 \mathbb E\bigl[\inf\{k\ge0:\|\nabla\psi(x_k)\|\le\varepsilon\}\bigr]
 \le\frac{16(\psi(x_0)-\psi_*)}{\alpha\varepsilon^2}\,.
 \label{eq:stopping-time-bound}
\end{equation}

To apply this framework, we take $\tilde g_k:=\hat g_h(x_k)$, computed at
each iteration with a new direction batch of either type in
\cref{def:spherical-directions} and new samples and comparison variables as
in~\eqref{eq:estimator}, and write
$K_\varepsilon:=\inf\{k\ge0:\|G(x_k)\|\le\varepsilon\}$.

\begin{theorem}
\label{thm:opt-comparison-stationarity}
Suppose that $f(\cdot;\xi)$ satisfies \cref{ass:function-smooth} for almost
every $\xi$, that $\nabla f=\mathbb E_\xi\nabla f(\cdot;\xi)$ and
$\mathbb E_\xi\|\nabla f(\cdot;\xi)\|^3\le\sigma_f^3$, and that
$f_*:=\inf f>-\infty$. For $\varepsilon>0$, let $h$ satisfy the condition
in \cref{cor:gradient-comparison-complexity} with $m_3$ replaced by
$\sigma_f^3$, and let
$0<\alpha\le\min\{\tfrac1{\gamma L_f},\tfrac{h^2s\varepsilon^2}{2\gamma L_fd^2}\}$.
Then the method $x_{k+1}=x_k-\alpha\hat g_h(x_k)$ reaches a point that is
$\varepsilon$-stationary with respect to comparisons, and
\[
 \mathbb E[K_\varepsilon]\le\frac{16\gamma(f(x_0)-f_*)}{\alpha\varepsilon^2}\,.
\]
\end{theorem}
\begin{proof}
We apply \citet[Theorem~4.18(1)]{berahas_cao_scheinberg_oracles_review},
recalled in~\eqref{eq:stopping-time-bound}, to $\psi:=\gamma f$ in place of
their objective $\phi$, with $\tilde g_k=\hat g_h(x_k)$ as the oracle output.
By \cref{lem:comparison-gradient}(ii), $\nabla\psi=G$, so the stopping time
in~\eqref{eq:stopping-time-bound} is $K_\varepsilon$. Their $L_1$ and
$\hat\phi$ become $L_\psi=\gamma L_f$ and $\psi_*=\gamma f_*$. By
\cref{thm:gradient-estimator} and the condition on $h$, which keeps the bias
of $\hat g_h$ at most $\tfrac\varepsilon2$, \eqref{eq:oracle-conditions} holds with $\epsilon_g^2=\tfrac{\varepsilon^2}{4}$ and
$\sigma_c^2=\tfrac{d^2}{4h^2s}$. The bound on $\alpha$ gives
$\alpha\le\tfrac1{L_\psi}$ and
\[
 \epsilon_g^2+\alpha L_\psi\sigma_c^2\le\frac{\varepsilon^2}4+\frac{\varepsilon^2}8=\frac{3\varepsilon^2}8\,,
\]
so~\eqref{eq:stopping-time-bound} gives
$\mathbb E[K_\varepsilon]\le\tfrac{16\gamma(f(x_0)-f_*)}{\alpha\varepsilon^2}$.
\end{proof}

\begin{remark}
\label[remark]{rem:opt-comparison-counts}
Each iteration uses $s$ comparisons. With the largest $\alpha$ in
\cref{thm:opt-comparison-stationarity}, the expected number of comparisons
satisfies
\[
 s\,\mathbb E[K_\varepsilon]
 \le\frac{16\gamma^2L_f(f(x_0)-f_*)}{\varepsilon^2}
 \max\Bigl\{s,\frac{2d^2}{h^2\varepsilon^2}\Bigr\}\,.
\]
For fixed problem parameters,
$K_{\mathcal D}=O(1)$, and small $\varepsilon$, the choices of $h$ in
\cref{rem:comparison-counts}, with $m_3$ replaced by $\sigma_f^3$, make this
$O(d^2\varepsilon^{-6})$, and $O(d^2\varepsilon^{-5})$ when $f(\cdot;\xi)$
also satisfies \cref{ass:higher-order-smooth}. The same framework also gives
rates under PL or linear-KL conditions
\citep{berahas_cao_scheinberg_oracles_review}.
\end{remark}

\section{Estimating the smoothed comparison gradient}
\label{sec:nonsmooth-estimation}

We now consider sample functions that are Lipschitz but not necessarily
differentiable. We replace \cref{ass:function-smooth,ass:higher-order-smooth}
by the following condition.

\begin{assumption}[Lipschitz sample functions]
\label{ass:nonsmooth-lipschitz}
There exists $L>0$ such that $|f(y;\xi)-f(z;\xi)|\le L\|y-z\|$ for all
$y,z\in\mathbb R^d$ and all $\xi$.
\end{assumption}

Without smoothness, the comparison gradient $G(x)$ is no longer a good
target. Comparisons near $x$ only reveal $G$ near $x$, and $G$ may change
abruptly. Indeed, the lower bound of
\cref{thm:comparison-complexity-lower-bounds}(i) grows as $L_f^2$, and its
hard family is $\tfrac{4\varepsilon}{\gamma}$-Lipschitz, so under
\cref{ass:nonsmooth-lipschitz} alone no bounded number of comparisons
suffices to estimate $G(x)$. Nonsmooth optimization faces the same
difficulty with $\nabla f$ and works instead with the gradient of a smoothed
function, such as the Gaussian smoothing
\[
 f_h(x):=\mathbb E_u[f(x+hu)]\,,\qquad u\sim N(0,\tfrac1dI_d)\,,
\]
with smoothing parameter $h>0$. \citet{NesterovSpokoiny2017} show that
$f_h$ is within $hL$ of $f$ and has a $\tfrac{dL}{h}$-Lipschitz gradient, and
use it in zeroth-order methods for nonsmooth convex and nonconvex problems.
\citet{DuchiBartlettWainwright2012} use randomized smoothing in stochastic
methods for nonsmooth convex problems. \citet{LinZhengJordan2022} relate
smoothing over a ball to the Goldstein subdifferential, and
\citet{ChaytiElBakkaliSaadi2026NoisyComparisons} use this relation to find
Goldstein stationary points of deterministic Lipschitz functions from
comparisons when $\Phi$ is known.

We do the same with comparisons. We average $G$ over the same Gaussian
neighborhood of $x$ and call
\[
 G_h(x):=\mathbb E_u[G(x+hu)]
\]
the \emph{smoothed comparison gradient}.
Like $G(x)$, it is defined by the responses of the oracle alone. For a unit
vector $v$ and the random point $x+hu$ near $x$, the oracle prefers
$x+hu+tv$ to $x+hu$ with probability $\tfrac12-t\langle G_h(x),v\rangle+o(t)$
as $t\to0$, on average over $u$. Thus $\|G_h(x)\|$ is the largest average
rate at which a small move from points near $x$ can increase the probability
of being preferred. The next lemma shows that $G_h(x)$ is $\gamma$ times the
gradient of $f_h$.

\begin{lemma}
\label[lemma]{lem:nonsmooth-identity}
Suppose that \cref{ass:nonsmooth-lipschitz} holds. Then, for any $h>0$,
$f_h$ is continuously differentiable and $G_h(x)=\gamma\nabla f_h(x)$ for all
$x\in\mathbb R^d$.
\end{lemma}
\begin{proof}
Since every $f(\cdot;\xi)$ is $L$-Lipschitz, $f_h$ is continuously
differentiable, and dominated convergence, with difference quotients bounded
by $L$, gives $\nabla f_h(x)=\mathbb E_{u,\xi}[\nabla f(x+hu;\xi)]$. For
almost every $u$, \cref{lem:comparison-gradient}(ii) gives
$G(x+hu)=\gamma\,\mathbb E_\xi\nabla f(x+hu;\xi)$, and taking expectations
over $u$ proves the claim.
\end{proof}

To estimate $G_h(x)$, we compare a random base point near $x$ with a
perturbation of it of a smaller size $\eta\le h$, as \citet{duchi2015optimal}
do with function values for nonsmooth convex problems.

\begin{estimatorbox}{Smoothed random-direction comparison estimator}
Let $x\in\mathbb R^d$ be an estimation point, $h>0$ a smoothing parameter,
$0<\eta\le h$ a sampling radius, and $1\le s\le d$ an integer. Draw
independent vectors $u_1,\ldots,u_s,w_1,\ldots,w_s\sim N(0,\tfrac1dI_d)$, and let
$y_i:=x+\sqrt{h^2-\eta^2}\,w_i$. For each $i$, draw a sample
$\xi_i\sim\mathbb P_\xi$ and a comparison variable $\zeta_i$, compare $y_i$
and $y_i+\eta u_i$ on the sample function $f(\cdot;\xi_i)$, and set
\begin{equation}
 \hat g_{h,\eta}(x):=\frac{d}{2\eta s}
 \sum_{i=1}^s C_{\zeta_i}(y_i,y_i+\eta u_i;\xi_i)u_i\,.
 \label{eq:nonsmooth-estimator}
\end{equation}
The samples and comparison variables are mutually independent and independent
of the $u_i$ and $w_i$.
\end{estimatorbox}
Since $(h^2-\eta^2)+\eta^2=h^2$, each $y_i+\eta u_i$ has the distribution of
$x+hu$, so the target $G_h(x)$ does not depend on $\eta$, which controls only
the estimation error.

\begin{theorem}
\label{thm:nonsmooth-estimator}
Suppose that \cref{ass:nonsmooth-lipschitz} holds. For any
$x\in\mathbb R^d$, $h>0$, $0<\eta\le h$, and $1\le s\le d$, the
estimator~\eqref{eq:nonsmooth-estimator} satisfies
\begin{equation*}
 \mathbb E\|\hat g_{h,\eta}(x)-\mathbb E\hat g_{h,\eta}(x)\|^2
 \le\frac{d^2}{4\eta^2s}\,,
 \qquad\text{and}\qquad
 \|\mathbb E\hat g_{h,\eta}(x)-G_h(x)\|\le\frac{M_\phi L^3}{2}\eta^2\,.
\end{equation*}
\end{theorem}
\begin{proof}
The summands in~\eqref{eq:nonsmooth-estimator} are independent with second
moment $\mathbb E\|u_i\|^2=1$, so the variance bound follows as in the proof
of \cref{lem:variance-estimator}.

For the bias, let $(y,u,\xi)$ have the distribution of $(y_i,u_i,\xi_i)$.
Since $y+\eta u$ has the distribution of $x+hu$,
\cref{lem:nonsmooth-identity} and its proof give
\begin{equation}
 G_h(x)=\gamma\,\mathbb E_{y,u,\xi}[\nabla f(y+\eta u;\xi)]\,.
 \label{eq:nonsmooth-smoothed-gradient}
\end{equation}
By~\eqref{eq:comparison-probability},
\[
 \mathbb E[C_\zeta(y,y+\eta u;\xi)\mid y,u,\xi]=2\Phi(f(y+\eta u;\xi)-f(y;\xi))\,.
\]
For $u\sim N(0,\tfrac1dI_d)$ and a Lipschitz
function $g$, Gaussian integration by parts gives
$d\,\mathbb E_u[g(u)u]=\mathbb E_u[\nabla g(u)]$. We apply it to
$g(u):=\Phi(f(y+\eta u;\xi)-f(y;\xi))$, whose gradient is
$\eta\Phi'(f(y+\eta u;\xi)-f(y;\xi))\nabla f(y+\eta u;\xi)$ by the chain
rule, and obtain
\begin{equation*}
 \mathbb E\hat g_{h,\eta}(x)
 =\frac{d}{\eta}\mathbb E_{y,u,\xi}\bigl[\Phi(f(y+\eta u;\xi)-f(y;\xi))u\bigr]
 =\mathbb E_{y,u,\xi}\bigl[\Phi'(f(y+\eta u;\xi)-f(y;\xi))\nabla f(y+\eta u;\xi)\bigr]\,.
\end{equation*}
Subtracting~\eqref{eq:nonsmooth-smoothed-gradient}, we obtain
\begin{align*}
 \|\mathbb E\hat g_{h,\eta}(x)-G_h(x)\|
 &=\bigl\|\mathbb E_{y,u,\xi}\bigl[\bigl(\Phi'(f(y+\eta u;\xi)-f(y;\xi))-\gamma\bigr)
   \nabla f(y+\eta u;\xi)\bigr]\bigr\|\\
 &\le\frac{M_\phi}{2}
   \mathbb E_{y,u,\xi}\bigl[|f(y+\eta u;\xi)-f(y;\xi)|^2\,
   \|\nabla f(y+\eta u;\xi)\|\bigr]
 \le\frac{M_\phi L^3\eta^2}{2}\,,
\end{align*}
where the first inequality uses $|\Phi'(t)-\gamma|\le\tfrac{M_\phi}2t^2$
from~\eqref{eq:Phi-bounds}, and the second uses
\cref{ass:nonsmooth-lipschitz} and $\mathbb E\|u\|^2=1$.
\end{proof}

As in \cref{cor:gradient-comparison-complexity}, averaging gives the number
of comparisons.

\begin{corollary}
\label[corollary]{cor:nonsmooth-complexity}
Suppose that \cref{ass:nonsmooth-lipschitz} holds. For any given $h>0$ and
$\varepsilon>0$, let $\hat g$ be the average of $T$ independent realizations
of~\eqref{eq:nonsmooth-estimator} with $s=1$. If $\eta\le h$,
$M_\phi L^3\eta^2\le\varepsilon$, and
$T\ge\tfrac{d^2}{2\eta^2\varepsilon^2}$, then
$\mathbb E\|\hat g-G_h(x)\|^2\le\varepsilon^2$.
\end{corollary}
\begin{proof}
By \cref{thm:nonsmooth-estimator} and independence,
\[
 \mathbb E\|\hat g-G_h(x)\|^2
 \le\frac{d^2}{4\eta^2T}+\frac{M_\phi^2L^6\eta^4}{4}
 \le\frac{\varepsilon^2}{2}+\frac{\varepsilon^2}{4}\le\varepsilon^2\,.\qedhere
\]
\end{proof}

\begin{remark}
\label[remark]{rem:nonsmooth-comparison-counts}
The largest $\eta$ allowed in \cref{cor:nonsmooth-complexity} is
$\min\{h,\sqrt{\tfrac{\varepsilon}{M_\phi L^3}}\}$. With this choice, $T$ can be
taken of order $d^2\varepsilon^{-2}\max\{h^{-2},M_\phi L^3\varepsilon^{-1}\}$,
which is $d^2\varepsilon^{-3}$ for fixed $h$ and small $\varepsilon$.
\end{remark}

\appendix
\section*{Appendix}
\phantomsection\label{app:proofs}

\begin{proof}[Proof of the bounds~\eqref{eq:Phi-bounds}]
The first two bounds follow
by integrating the bound on $\Phi'''$ and using $\Phi''(0)=0$ and
$\Phi'(0)=\gamma$. For the third, by Taylor's theorem,
$\Phi(t+\delta)-\Phi(t-\delta)=2\delta\Phi'(t)+r$ with
$|r|\le\tfrac{M_\phi\delta^3}{3}$. Because $|\Phi|\le\tfrac12$, the choice
$\delta=M_\phi^{-1/3}$ yields
\[
 |\Phi'(t)|\le\frac{1}{2\delta}+\frac{M_\phi\delta^2}{6}
 =\frac23M_\phi^{1/3}\,.\qedhere
\]
\end{proof}

\begin{proof}[Proof of \cref{lem:hard-family-smoothness}]
We write $D^kg(z)$ for the $k$th derivative of a map $g$ at $z$, which is a
$k$-linear map, and extend $\|\cdot\|_{\mathrm{op}}$ to a $k$-linear map $T$
by $\|T\|_{\mathrm{op}}:=\sup\{\|T[u_1,\ldots,u_k]\|:\|u_1\|,\ldots,\|u_k\|\le1\}$.
Since $\Psi_\rho(z)=\rho\Psi_1(\tfrac z\rho)$, the Cauchy--Schwarz
inequality gives, for $k=2,3$,
\begin{equation*}
\left\lVert D^k f_\theta(z)\right\rVert_{\mathrm{op}}
\le\left\lVert\theta\right\rVert\rho^{1-k}
\left\lVert D^k\Psi_1\Bigl(\frac{z-x}{\rho}\Bigr)\right\rVert_{\mathrm{op}}\,.
\end{equation*}
Integrating along segments, it therefore suffices to show
$\lVert D^2\Psi_1(z)\rVert_{\mathrm{op}}<3$ and
$\lVert D^3\Psi_1(z)\rVert_{\mathrm{op}}<13$ for all $z$.
Let $r:=\lVert z\rVert$ and let $u,v,w$ be unit vectors. Direct
differentiation gives
\begin{equation*}
D^2\Psi_1(z)[u,v]
=-\frac{2(\langle z,u\rangle v+\langle z,v\rangle u+\langle u,v\rangle z)}{(1+r^2)^2}
+\frac{8\langle z,u\rangle\langle z,v\rangle z}{(1+r^2)^3}\,,
\end{equation*}
whose two vector factors have norms at most $3r$ and $r^3$. Differentiating
in direction $w$ produces three terms with coefficient
$-\tfrac{2}{(1+r^2)^2}$, six with $\tfrac{8}{(1+r^2)^3}$, and one with
$-\tfrac{48}{(1+r^2)^4}$, whose vector factors have norms at most $1$,
$r^2$, and $r^4$. Hence
\begin{align*}
\lVert D^2\Psi_1(z)\rVert_{\mathrm{op}}
&\le\frac{6r}{(1+r^2)^2}+\frac{8r^3}{(1+r^2)^3}<3\,,
\qquad\text{and}\\
\lVert D^3\Psi_1(z)\rVert_{\mathrm{op}}
&\le\frac{6}{(1+r^2)^2}+\frac{48r^2}{(1+r^2)^3}+\frac{48r^4}{(1+r^2)^4}<13\,,
\end{align*}
where the strict inequalities follow by maximizing over $r\ge0$.
\end{proof}

\bibliographystyle{abbrvnat}
\bibliography{reference}

\end{document}